\documentclass[preprint,authoryear,12pt]{elsarticle}
\usepackage{amsfonts,amssymb,amsmath,amsopn}
\usepackage{rotating}
\usepackage{bigints}
\usepackage{booktabs}
\usepackage{subfigure}
\usepackage{url}
\usepackage{amssymb}
\usepackage{color}
\usepackage{multirow}
\usepackage{diagbox}
\usepackage{enumerate}
\usepackage{enumitem}
\usepackage{epstopdf}
\usepackage{float}
\usepackage{graphicx}
\usepackage{caption}
\usepackage[colorlinks=true]{hyperref}
\newcommand{\ignore}[1]{}
\newtheorem{asu}{{\sc Assumption}}

\newtheorem{pro}{{\sc Proposition}}
\newtheorem{thm}{Theorem}

\newdefinition{rmk}{Remark}

\usepackage[margin=1in]{geometry}
\newproof{pf}{Proof}
\newproof{pot}{Proof of Theorem \ref{thm2}}

\usepackage{amssymb}
\usepackage{booktabs}
\usepackage{bm}

\usepackage{xr}
\graphicspath{{figure/}} 

\newcommand{\E}{\mathbb{E}}

\newcommand{\red}{}
\newcommand{\tr}{\text{tr}}
\newcommand{\T}{\mathsf{T}}

\allowdisplaybreaks

\begin{document}
	\begin{frontmatter}
		\title{
Spectral analysis of high-dimensional spot volatility matrix with applications
		}
		\author[SHFE]{Qiang LIU\thanks{cor1}}
		\author[JNU]{Yiming LIU}
		\author[UM]{Zhi LIU}
		\author[NUS]{Wang ZHOU}
		\cortext[cor1]{Corresponding author. Email: liuqiang@mail.shufe.edu.cn. \ The authors are listed in alphabetic order.}
		\address[SHFE]{School of Statistics and Data Science, Shanghai University of Finance and Economics}
		\address[JNU]{School of Economics, Jinan University}
		\address[UM]{Department of Mathematics, University of Macau}
		\address[NUS]{Department of Statistics and Data Science, National University of Singapore}
		\begin{abstract}
In random matrix theory, the spectral distribution of the covariance matrix has been well studied under the large dimensional asymptotic regime when the dimensionality and the sample size tend to infinity at the same rate. However, most existing theories are built upon the assumption of independent and identically distributed samples, which may be violated in practice. For example, the observational data of continuous-time processes at discrete time points, namely, the high-frequency data. In this paper, we extend the classical spectral analysis for the covariance matrix in large dimensional random matrix to the spot volatility matrix by using the high-frequency data. We establish the first-order limiting spectral distribution and obtain a second-order result, that is, the central limit theorem for linear spectral statistics. Moreover, we apply the results to design some feasible tests for the spot volatility matrix, including the identity and sphericity tests. Simulation studies justify the finite sample performance of the test statistics and verify our established theory. 
		\end{abstract}
	\end{frontmatter}
	
\section{Introduction}

During the past decades, it has been a hot topic for investigating the asymptotic spectral analysis of sample covariance matrices when the data dimension increases to infinity at the same rate as the sample size increases (Such a scenario is called large dimensional setting). It is shown that the theory differs from the classical limiting theory in the fixed dimension. To be specific, the limiting spectral distributions of the population covariance matrix and the sample covariance matrix will be governed by the so-called Mar$\check{\text{c}}$enko-Pastur equation (\cite{S1995}), instead of the empirical central limit theorem at the fixed dimension. Interested readers can refer to \cite{B2010} for a comprehensive introduction of the large dimensional random matrix theory. However, most existing works assume that the samples are independent and identically distributed, which is not satisfied in many situations, especially for time series data. In this paper, the data type to our concern is the high-frequency data, which is neither independent nor identically distributed (Such a point will be elaborated in Section \ref{sec:setup}).  

The computational technology has developed fast and has been widely applied in the financial market in recent years, making the high-frequency trading strategies possible and generating massive high-frequency data. As a result, developing statistical models and econometric methods for the high-frequency data has been experiencing exponential growth, both by practicers and researchers. The key objectives that attract most attention from statisticians and econometricians are the volatility of a single asset and the volatility matrix for the multivariate case (both of these two quantities have the integrated version and the spot version), since they play pivotal roles in many areas of financial economics, including asset and derivative pricing, portfolio allocation, risk management, hedging, and etc. 
Various methods of estimation and statistical inference regarding to these mentioned quantities have been proposed, and a bulk of related references can be found in \cite{AJ2014}. We notice that most of the existing literature focuses on one-dimensional case or fixed dimension, while the consideration of high dimension allowing the dimensionality tend to infinity is rather rare in comparison. It is well known that the fact of high dimension is a common feature for vast datasets in this big data era, thus the extension of classical theory under fixed dimension to this scenario becomes an important research direction in statistics. For the high-frequency data, when both the number of assets and the sample size of the transactional price data of the assets go to infinity, the estimation of integrated and spot volatility matrix has been considered in \cite{WZ2010}, \cite{TWC2013}, \cite{TWZ2013}, \cite{KKLW2018}, \cite{KLW2018}, \cite{K2018}, \cite{DLX2019}, \cite{BMN2019}, \cite{BDLW2022}, \cite{LLZ2024}, and references therein. 
Instead of directly estimating the integrated and spot volatility matrix, their spectral analysis also has wide application in multivariate hypothesis testing, principal component analysis, and factor analysis, and is becoming a hot topic recently (See, e.g. \cite{AX2019}, \cite{CMZ2020}, \cite{AX2017}, \cite{K2023},  \cite{KLLL2023}.). For example, as suggested in \cite{CMZ2020}, investors should be better off investing in a statistically estimated principal component rather than an index fund. 
Under the large dimensional regime, namely the data dimension and the sample size tend to infinity at the same rate, \cite{ZL2011} studied the limiting spectral distribution of the realized volatility matrix estimator. They found that its limiting behavior is greatly affected by the time variability of the volatility matrix process, and the empirical result for independently and identically distributed samples established in large dimensional random matrix can not be directly used. Moreover, they proposed a time-variation adjusted realized volatility matrix estimator, whose limiting spectral distribution depends solely on that of the targeting integrated volatility matrix via a Mar$\check{\text{c}}$enko-Pastur equation, making the inference of the latter one to be possible. 
The same problem was also considered by \cite{HP2014}, who obtained the explicit form of the moments of the limiting spectral distribution of the realized volatility matrix estimator by using the tools from graph theory.
Based on \cite{ZL2011}, \cite{XZ2018} considered the further presence of market microstructure noise in the high-frequency data and dealt with its effect by pre-averaging technique. 
Furthermore, for the specific volatility matrix process considered in \cite{ZL2011}, \cite{YZC2021} established the central limit theorem for the linear spectral statistics of the time-variation adjusted realized covariance matrix estimator proposed in \cite{ZL2011}.  Based on this, they proposed several test statistics for the identity and sphericity hypotheses of the integrated volatility matrix.

In this paper, our focus is the spot volatility matrix, which can quantify the co-variation pattern of the price processes of multiple assets at any given time second. Unlike the integrated volatility matrix, which is defined over a fixed time interval, say a day, the spot volatility matrix is a more general and flexible object, leading itself to wider range of applications.
\red{For example, in practice, most financial and macroeconomic time series exhibit time-varying volatility, and the expected volatility matrix in the near future can be very different from the average of the volatility matrix over a long time horizon. To reduce the bias of the forecasted volatility matrix, which plays a pivotal role in mean-variance efficient portfolio theory, we need to shorten the learning period to better capture the dynamics of the spot volatility matrix (See, e.g. \cite{FLY2012}, \cite{JLT2024}, and references therein).} 
The first question we are trying to figure out is the relationship between the limiting spectral distribution of the spot volatility matrix and the one of the realized spot volatility matrix estimator, under the infill setting of high-frequency data, namely, the length of time interval between two consecutive observations decreases to 0. Meanwhile, we allow the data dimension tends to infinity at the same rate as the sample size increases. For identically and independently distributed data, it is well known that the limiting spectral distribution of the population covariance matrix is linked with one of the sample covariance matrix by a Mar$\check{\text{c}}$enko–Pastur equation through Stieltjes transform when the date dimension and the sample size tend to infinity at the same rate (\cite{S1995}). 
Our analysis shows that such a result remains valid for the spot volatility matrix and its estimator under some regular conditions of smoothness for the volatility matrix process. 
Furthermore, we also investigate the central limit theorem for the linear spectral statistics of the spot volatility matrix, which serves as a second-order limiting result and has wide applicability in multivariate statistical inference problems. 
As for possible applications of our established theoretical results, we consider two hypothesis testing problems on if the spot volatility matrix is equal or proportional to a given matrix, which are called identity test and sphericity test respectively. To this end, we propose several test statistics in the large dimensional regime and demonstrate their finite sample performances via simulation studies.  

The rest of the article is organized as follows. In Section \ref{sec:setup}, we introduce the framework of high-frequency data, the related assumptions and present the realized spot volatility matrix estimator. In Section \ref{sec:asy}, we present the asymptotic results regarding to the empirical spectral distribution and the linear spectral statistics of the realized spot volatility matrix estimator. For the application of identity test and sphericity test for the spot volatility matrix, we propose several test statistics and present their asymptotic properties in Section \ref{sec:app}. Section \ref{sec:conc} concludes the paper. In a separate supplementary material, we conduct simulation studies to verify our established theory and demonstrate the finite sample performance of the proposed test statistics. 

\textbf{Notations:} At the end of this introduction part, we introduce some notation and definitions used throughout the paper. For any matrix $\Sigma$, we use $\Sigma^{(i,j)}$, $\Sigma^{\T}$,  $\tr (\Sigma)$ to denote its $(i,j)$-th element, transpose, and trace, respectively. We use $\lambda_i^{\Sigma}$, $\lambda_{\min}^{\Sigma}$, and $\lambda_{\max}^{\Sigma}$ to denote the $i$-th (in descending order), the minimum, and the maximum eigenvalues of a $p-$dimensional squared matrix $\Sigma$, respectively. Furthermore, the notation $\lambda(\Sigma) = (\lambda_1^{\Sigma},...,\lambda_p^{\Sigma})^{\T}$ means the function extracting the eigenvalues of $\Sigma$ as a $p$-dimensional vector, in a non-increasing order.
\red{Meanwhile, we let $\|\Sigma\|$ and $\|\Sigma\|_{F}$ be the spectral norm and Frobenius norm of $\Sigma$, namely $\|\Sigma\| := \sqrt{\lambda_{\max}^{\Sigma^{\T} \Sigma}}$ and $\|\Sigma\|_{F}:= \sqrt{\sum_{i=1}^{p} (\lambda_{i}^{\Sigma})^2}$. 
}
For any vector $a$, we use $a^{(j)}$ to denote its $j-$th entry. Specifically, we define $\mathbf{0}_{p}, \mathbf{1}_{p}$ as an $p-$dimensional vector with each entry being 0 and 1, respectively, and $\mathbb{I}_p$ as a $p \times p$ identity matrix. The formula $A \preceq B$ means that $B-A$ is strictly positive definite. 
For any interval $I \subset [0, \infty)$ and any metric space $S$, $D(I; S)$ stands for the space of c$\grave{\text{a}}$dl$\grave{\text{a}}$g functions from $I$ to $S$. We use the notation $\longrightarrow^d, \longrightarrow^{p}$ as the convergence in distribution and convergence in probability respectively, and $N(\mu,\sigma^2)$ as the normal distribution with mean $\mu$ and variance $\sigma^2$. \red{The notation $\lfloor x\rfloor$ stands for $\max \{k \in \mathbb{Z}: k\leq x\}$.}


\section{Setup and assumptions}\label{sec:setup}

Without the loss of generality, we define all the processes involved in the time interval $[0,1]$, where the time unit may be one day, one month, or one year. On the filtered probability space $(\Omega, \mathcal{F}, (\mathcal{F}_t)_{0\leq t\leq 1}, \mathbf{P})$, we denote the underlying efficient log-price process of the assets as $\{X_t\}_{0\leq t\leq 1}$. Moreover, we assume $X \in \mathbb{R}^{p}$, with $p$ being the total number of assets under consideration. 
According to the fundamental theorem of asset pricing (See, e.g. \cite{DS1994}), in a frictionless financial market with no-arbitrage opportunity, the price process of an asset is necessarily to be a semimartingale. 
In accordance with this, we assume that $X$ follows a continuous Brownian semimartingale and can be written as 
\begin{align}\label{model}
X_t = X_0 + \int_{0}^{t} b_sds +  \int_{0}^{t} \sigma_{s}dB_s, \quad t\in[0,1].
\end{align}
In the above model, $b \in \mathbb{R}^{p}$ is the drift term being progressively measurable and locally bounded, $B$ is a $q$-dimensional standard Brownian motion and its dimension $q$ usually means the number of risk factors involved in financial asset prices $X$ or the dimension of the state space of a continuous-time factor model (See, e.g. \cite{JM2013rank}, \cite{AX2017} and many others), $\sigma$ is a $p\times q$ matrix-valued stochastic process generating adapted and c$\grave{\text{a}}$dl$\grave{\text{a}}$g path almost surely.
In \eqref{model}, we allow $\sigma$ and $B$ to be mutually dependent with any general structure, thus depicting the leverage effect in finance.
We define the spot volatility matrix at any given time point $t\in[0,1]$ as $c_t  = \sigma_t \sigma_t^{\mathsf{T}}$ and its integral version, the integrated volatility matrix, as $\int_{0}^{1} c_tdt$\footnote{Without causing confusion, we will call both the processes $\sigma$ and $c$ as volatility process, from here and after.}.
Both measurements quantify the joint variational strength among different assets and play pivotal roles in financial economics.
In this paper, we are interested in exploring the local behavior of the volatility matrix process, thus we mainly focus on the spot volatility matrix.
More importantly, our framework allows the dimensionality $p$ to diverge with the sample size, in contrast to most existing high-frequency studies where $p$ is treated as fixed.
Besides, the factor number $q$ can be either fixed or divergent at a proper rate as $p\rightarrow \infty$.
For notational simplicity, we omit the dependence of the dimension parameters $p$ and $q$ for $\sigma_t, c_t$. We assume that our model \eqref{model} satisfies the following assumption: 

\begin{asu}\label{asu:dgp}
\red{
For $1\leq i\leq p$ and $1\leq j\leq q$, the components of $b_t^{(i)}$ and $\sigma_t^{(i,j)}$ are adapted locally bounded processes with continuous sample path. 
Moreover, for almost all path of $\sigma_t^{(i,j)}$ within $t\in[0,1]$, it holds that 
	\begin{equation}\label{asu:dgp:cont}
			 \sup_{1\leq i\leq p, 1 \leq j \leq q} \E\left[\left|\sigma^{(i,j)}_{t+\epsilon}-\sigma^{(i,j)}_{t} \right|^{\red{2}} \right]  \leq B(t,\epsilon) |\epsilon|^{\red2\gamma} + o(|\epsilon|^{\red2\gamma}), \ \epsilon \rightarrow 0,
	\end{equation}
for {\red{some $\gamma>0$}}, where $B(t,\epsilon)$ is a positive random function and continuous with respect to $t$. 
And there exists a strictly positive definite matrix $\underline{c}$ such that $\underline{c} \preceq c_t$ holds uniformly over the time interval $[0,1]$, almost surely.
}
\end{asu}

\begin{rmk}
Assumption \ref{asu:dgp} essentially imposes some local boundedness and continuity conditions for the entries of the drift process $b$ and the volatility matrix process $c$. According to the localization procedure described in Section 4.4.1 of \cite{JP2012}, these local conditions can be equivalently turned to corresponding boundedness and continuity conditions, which will be directly used in our proof. 
Condition \eqref{asu:dgp:cont} requires that each entry of $\sigma$ to be as smooth enough, such condition will be used to approximate the spot volatility matrix at any given time $t\in[0,1]$ by the ones closed to time $t$, and the approximation error will depend on the smoothness parameter $\gamma$.  
It is obvious that as $\gamma$ increases, the generated paths of $\sigma^{(i,j)}$ with $1 \leq i \leq p, 1 \leq j\leq q$ are relative smoother, almost surely. In this sense, the case of $\gamma = \infty$ is the most restrictive one, and only the constant volatility model can meet the requirement, while such a continuity condition plays no extra role if we set $\gamma = 0$.
When $0<\gamma \leq \frac{1}{2}$, a specific model of continuous Brownian semimartingale like \eqref{model}, which includes the famous Heston model, or general It$\hat{\text{o}}$ semimartingale incorporating a jump part is often used for the volatility process (See, e.g. \cite{JT2014, LL2023} and many others), due to its wide application in finance. \red{Some literature in finance report that the volatility path can be much rough than a semimartingale, showing a Holder continuous path of order about $0.1$, which is the value of $\gamma$. This is called the rough volatility, a very hot topic in recent ten years (See, e,g. \cite{GJR2018}, \cite{CHLR2024}, \cite{CT2025}, and references therein).}
We also note that a smoothness condition similar to \eqref{asu:dgp:cont} is adopted in \cite{BDLW2022} for the estimation of spot volatility matrix and in \cite{LLL2018}, \cite{LL2022} for the univariate case of spot volatility.
\end{rmk}
  
In practice, the continuous sample path of $X_t$ over $[0,1]$ is unobservable, and it can only be observed at some discrete time points. 
Throughout this paper, we consider the equidistant time grids scattered over the interval $[0,1]$ and write them as $\{i/n: i=0,1,...,n\}$, where $n$ is the total number of observations. 
We see that as $n$ increases, the time length between two consecutive observations shrinks, resulting in the so-called high-frequency data. Our theoretical results will be established under the infill asymptotic setting of $n\rightarrow \infty$. 
With the definition $\Delta_i^n X = X_{i/n} - X_{(i-1)/n}$ for $i=1,..., n$, it is well known that the integrated volatility matrix $\int_{0}^{1} c_td_t$ can be estimated by the realized volatility matrix estimator $\widehat{IV}^n$:
\begin{align}\label{intvolest}
	\widehat{IV}^n =  \sum_{i=1}^{n} \left(\Delta_i^n X\right) \left(\Delta_i^n X\right)^{\mathsf{T}}.
\end{align}
And its local version of spot volatility matrix estimator at any given time $t \in [0,1]$ can be written as 
\begin{align}\label{spotvolest}
\widehat{c_{t}}^{n}  =\frac{n}{k_n}  \sum_{i= \lfloor tn \rfloor +1}^{ \lfloor tn \rfloor +k_n} \left(\Delta_i^n X\right) \left(\Delta_i^n X\right)^{\mathsf{T}},
\end{align}
where $k_n$ is a sequence of integers tending to infinity as $n \rightarrow \infty$, and meanwhile $k_n/n \rightarrow 0$.
A complete introduction and theoretical analyses of these two estimators can be found in \cite{AJ2014}, when the dimension $p$ is fixed as $n\rightarrow \infty$. 

\red{
It is worthy noting that under \eqref{model}, $\{\Delta_i^n X: i=1,...,n\}$ are neither independent nor identically distributed. To see this, if the stochastic processes $b$ and $\sigma$ are independent of the Brownian motion $B$, since $\Delta_i^n X = \int_{(i-1)/n}^{i/n} b_sds +  \int_{(i-1)/n}^{i/n} \sigma_sdB_s $, we have
\begin{align*}
 \Delta_i^n X  \sim MN\left(\int_{(i-1)/n}^{i/n} b_sds, \int_{(i-1)/n}^{i/n} c_sds\right),
\end{align*}
where $MN$ stands for mixed normal distribution.
Since $b$ and $c$ can be both autocorrelated and time-varying, $\Delta_i^n X$'s are not independent.
In random matrix theory, the spectral properties of the sample covariance matrix have been established by using $i.i.d.$ data under the large dimensional regime (\cite{B2010}).
In this paper, we are interested in extending the theory to the study of spot volatility matrix based on similar technique called Stieltjes transform. Another possible alternative method to solve this problem is to treat the spot volatility matrix estimator \eqref{spotvolest} as matrix-valued realization of $X$ driven by the stochastic differential equation \eqref{model} and investigate its spectral properties as in \cite{MP2022}, \cite{SYY2020}, \cite{SYY2021}, \cite{SYY2022} and references therein. The details remain to be explored.  
}

\section{Asymptotic results}\label{sec:asy}

In this section, we present the results of empirical spectral distribution (ESD) and linear spectral statistics (LSS) for the realized spot volatility matrix estimator $\widehat{c_{t}}^{n}$. The definitions of ESD and LSS will be given later. In the sequel, we first review the empirical results of ESD and LSS for covariance matrix under the large dimensional regime, which are followed by the extension to spot volatility matrix. From here and after, for $i=1,..., n$, we define $\Delta_i^n S= S_{i/n} - S_{(i-1)/n}$ for a general process $S$. We use the notation $C$ for a general positive constant, which may take different values from line to line.
 
\subsection{Empirical spectral distribution (ESD)}\label{sec:esd}
Let $M$ be a $p\times p$ random, symmetric and nonnegative definite matrix, whose eigenvalues are real and denoted as $\{\lambda_i^{M}: i=1,..., p\}$. Then, the empirical spectral distribution (ESD) of $M$ is defined as 
\begin{align}
F^{M}(x) := \frac{1}{p} \sum_{i=1}^{p} \mathbf{I}_{\{ \lambda_i^{M} \leq x\}}, \ x\in \mathbb{R}, 
\end{align}
where $\mathbf{I}_{\{\cdot\}}$ is the indicator function. In large dimensional random matrix theory, the relationship between the limiting spectral distribution of the population covariance matrix and the ESD of the sample covariance matrix is given by the Mar$\check{\text{c}}$enko–Pastur equation through the Stieltjes transform. 

\begin{pro}[Theorem 1.1 of \cite{S1995}]\label{pro1} 
Assume on a probability space, the following conditions are satisfied:
	\begin{enumerate}[label=(\roman*)]
\item for $p=1,2,...$ and for $1 \leq l\leq n$, $\mathbf{Z}_{l}^{(p)} = (Z_l^{(p,j)})_{1\leq j\leq p}$ with $Z_l^{(p,j)}$ $i.i.d.$ with mean 0 and variance 1;
\item as $n \rightarrow \infty$, it holds that $p \rightarrow \infty$ with $p/n \rightarrow y>0$;
\item $\Sigma_p$ is a (possibly random) nonnegative definite $p \times p$ matrix such that its ESD $F^{\Sigma_p}$ converges almost surely to a probability distribution $H$ on $[0, \infty)$ as $p \rightarrow \infty$;
\item $\Sigma_p$ and $\mathbf{Z}_{l}^{(p)}$'s are independent.
	\end{enumerate}
Let $\Sigma_p^{1/2}$ be the (nonnegative) square root matrix of $\Sigma_p$ and $S_p := 1/n \times \sum_{l=1}^{n} \Sigma_p^{1/2} \mathbf{Z}_l^{(p)} (\mathbf{Z}_l^{(p)})^{\T} \Sigma_p^{1/2}$. Then, almost surely, the ESD of $S_p$ converges weakly to a probability distribution $F^{y,H}$, which is determined by $H$ in that its Stieltjes transform
\begin{align}\label{equ:st}
m(z) := \int_{\lambda \in \mathbb{R}} \frac{1}{\lambda - z} dF^{y,H}(\lambda), \qquad z\in \mathbb{C}_{+} := \{z\in\mathbb{C}: Im(z) > 0\} 
\end{align}
is the only solution to the following Mar$\check{\text{c}}$enko-Pastur equation 
\begin{align}\label{equ:esd}
m(z) =  \int_{\tau \in \mathbb{R}} \frac{1}{\tau (1-y(1+zm(z)) )-z}dH(\tau), \quad z\in\mathbb{C}_{+},
\end{align}
in the set $\{m\in \mathbb{C}: -(1-y)/z + cm \in \mathbb{C}_{+}\}$.
\end{pro}

Specifically, for the special case of $\Sigma_p = \sigma^2 \mathbb{I}_p$ with some positive random variable $\sigma^2$, the limiting spectral distribution $F^{y, H}$ in the above proposition has an explicit expression called Mar$\check{\text{c}}$enko-Pastur law. Conditional on the ratio index $y$ and the scale index $\sigma^2$, the probability density function of $F^{y,H}$ is given by
\begin{align*}
p(x) = \frac{1}{2\pi \sigma^2 xy} \sqrt{(b-x)(x-a)}, \quad \text{if} \quad a\leq x\leq b,
\end{align*}
where $a=\sigma^2(1-\sqrt{y})^2, b=\sigma^2(1+\sqrt{y})^2$, and a point mass $1-1/y$ at the origin if $y>1$.
Let $\underline{m}(z)$ be the Stieltjes transform of $\underline{F}^{y,H} = y F^{y,H} + (1-y) \mathbf{1}_{[0, \infty)}$, which is the limiting ESD of $\underline{S}_p := 1/n \times \sum_{l=1}^{n} (\mathbf{Z}_l^{(p)})^{\T} \Sigma_p \mathbf{Z}_l^{(p)}$. Then, the equation \eqref{equ:esd} can be equivalently written as
\begin{align*}
z = -\frac{1}{\underline{m}(z)} + y \int \frac{t}{1+t\underline{m}(z)}dH(t), \quad z\in\mathbb{C}_{+},
\end{align*}
which is called Silverstein's equation, proposed in \cite{S1995}.

As for the high-frequency setting considered  in this paper, if there is no drift term and the volatility process is constant  in \eqref{model}, namely
\begin{align}\label{constvol}
 \quad  b_t\equiv 0 \ \text{and} \ \sigma_t \equiv \sigma_0, \ \text{for} \ 0 \leq t\leq 1,
\end{align}
then, for $i=1,..., n$, $\sqrt{n} \Delta_i^n X$ are independent and identically distributed with the distribution $(\sigma_0 \sigma_0^{\T})^{1/2} \cdot N(\mathbf{0}_{p},\mathbb{I}_p)$. As a direct application of Proposition \ref{pro1}, the limiting spectral distributions of the integrated volatility matrix $\int_{0}^{1} c_td_t$ and the realized volatility matrix estimator \eqref{intvolest} are determined by \eqref{equ:esd} if $p/n \rightarrow y$ as $n \rightarrow \infty$.
For the general stochastic volatility model as in \eqref{model}, since the $i.i.d.$ condition for the increments $\{\Delta_i^n X:i=1,...,n\}$ is usually violated, as a consequence, the theoretical results established in Proposition \ref{pro1} cannot be directly applied.
In fact, in Proposition 3 of \cite{ZL2011}, with the consideration of the special case
\begin{align}
\quad  b_t\equiv 0 \ \text{and} \ \sigma_t = \gamma_t \mathbb{I}_p, \ \text{for} \ 0 \leq t\leq 1,
\end{align}
where $\gamma_t$ is a time-varying and nonrandom scalar, it is shown that the ESD of the realized volatility estimator \eqref{intvolest} does not converge to Mar$\check{\text{c}}$enko–Pastur law.
Based on this special example, we see that, under the high-frequency setting, the relationship between the limiting spectral distributions of the integrated volatility matrix $\int_{0}^{1} c_td_t$ and the realized volatility matrix estimator \eqref{intvolest} is not governed by \eqref{equ:esd} anymore. \cite{ZL2011} pointed out that this is caused by the time-variability of the volatility process and they did some detailed analyses on how the time-variability in the volatility process affects the limiting spectral distribution of the realized volatility matrix estimator. Moreover, for model \eqref{model} satisfying the following class $\mathcal{C}$:
\begin{align}\label{vol:classC}
\sigma_t = \gamma_t \Lambda_p, \ \text{for} \ 0 \leq t\leq 1,
\end{align}
where, $\gamma_t$ is a one-dimensional random process belongs to $ D([0,1]; \mathbb{R})$, $\Lambda_p$ is a $p\times p$ nonrandom matrix satisfying $\text{tr}(\Lambda_p \Lambda_p^{\T}) = p$, they proposed an alternative estimator of the integrated volatility matrix that can be used to infer the limiting spectral distribution of the integrated volatility matrix. For general stochastic volatility model \eqref{model} without any restriction, the relationship between the limiting spectral distributions of the integrated volatility matrix and its realized estimator remains unknown.

For the spot volatility matrix, similar to the above analyses, if the condition \eqref{constvol} is satisfied, Proposition \ref{pro1} implies that the limiting spectral distributions of the spot volatility matrix $c_t$ and the realized spot volatility matrix estimator \eqref{spotvolest} are determined by \eqref{equ:esd}, if the dimension $p$ and the sample size $k_n$ satisfy $p/k_n \rightarrow y$ as $k_n \rightarrow \infty$. Our first finding is that, with Assumption \ref{asu:dgp} and under some suitable asymptotic conditions, such a result remains true for the general stochastic volatility model \eqref{model}. Specifically, we have

\begin{thm}\label{spot_esd}
Under Assumption \ref{asu:dgp} and the following conditions:
	\begin{enumerate}[label=(\roman*)]
		\item as $n \rightarrow \infty$, $p \rightarrow \infty$ and meanwhile $q/n \rightarrow 0$, $qp^{\gamma}/n^{\gamma} \rightarrow 0$ and $p/k_n \rightarrow \bar{p}>0$;
		\item for the nonnegative definite $p \times p$ spot volatility matrix (possibly random) $c_t$ at any fixed time $t\in[0,1]$, its ESD $F^{c_t}$ converges almost surely to a probability distribution $H_t$ on $[0, \infty)$ as $p \rightarrow \infty$;
	\end{enumerate}
Then, condition on $\mathcal{F}_t$ and in probability, the ESD of $\widehat{c_{t}}^{n}$ converges weakly to a probability distribution $F^{\bar{p},H_t}$, which is determined by $H_t$ in that its Stieltjes transform $m(z)$ as defined in \eqref{equ:st}
	is the only solution to Mar$\check{\text{c}}$enko-Pastur equation \eqref{equ:esd} with $H$ being replaced by $H_t$.  
\end{thm}
\textit{Proof.} \ignore{We first note that we can assume that the two processes $\sigma$ and $B$ are mutually independent without loss of generality if only condition $(\romannumeral3)$ is satisfied. This is because, if we let $B'$ be the result of replacing $\{B^{(j)}: j\in \mathcal{I}_p\}$ with independent standard Brownian motions that are also independent of $\sigma$ and define
\begin{align*}
	\widehat{c_{t}}'^{n}  =\frac{n}{k_n}  \sum_{i=\lfloor tn \rfloor+1}^{\lfloor tn \rfloor+k_n} \left(\int_{(i-1)\Delta_n}^{i\Delta_n} b_sds + \int_{(i-1)\Delta_n}^{i\Delta_n} \sigma_sdB'_s\right) \left(\int_{(i-1)\Delta_n}^{i\Delta_n} b_sds + \int_{(i-1)\Delta_n}^{i\Delta_n} \sigma_sdB'_s\right)^{\mathsf{T}}.
\end{align*}
Recall that $\widehat{c_{t}}^{n}$ is defined as in \eqref{spotvolest}, then under condition $(\romannumeral3)$ and according to Lemma 2.2 in \cite{B1999}, we can obtain
\begin{align*}
	\|F^{\widehat{c_{t}}^{n}} - F^{\widehat{c_{t}}'^{n}} \| \leq \frac{rank(\widehat{c_{t}}^{n}-\widehat{c_{t}}'^{n})}{p} \leq \frac{2\eta_p}{p} \rightarrow 0, \quad \text{as} \ p\rightarrow \infty. 
\end{align*}
Namely, $\widehat{c_{t}}^{n}$ and $\widehat{c_{t}}'^{n}$ have the same limiting spectral distribution. 
And in the preceding proof steps, we shall directly assume that $\sigma$ is independent of $B$.
}
We define 
\begin{align}\label{def:sigtilde}
	\widetilde{c_{t}}^{n}  =\frac{n}{k_n}  \sum_{i=\lfloor tn \rfloor+1}^{\lfloor tn \rfloor+k_n} \left(\sigma_{t} \Delta_i^n B\right) \left(\sigma_{t} \Delta_i^n B\right)^{\mathsf{T}},
\end{align}
and, for $j=1,...,k_n$,
\begin{align}
	u_j =\frac{\int_{(\lfloor tn \rfloor+j-1)/n}^{(\lfloor tn \rfloor+j)/n}b_sds + \int_{(\lfloor tn \rfloor+j-1)/n}^{(\lfloor tn \rfloor+j)/n} (\sigma_s - \sigma_{t})dB_s }{\sqrt{k_n/n}}, \quad v_j = \frac{\sigma_{t} \Delta_{\lfloor tn \rfloor+j}^nB}{\sqrt{k_n/n}},
\end{align}
and 
\begin{align}\label{def:UV}
	U=(u_1+v_1,...,u_{k_n}+v_{k_n} ), \qquad V= (v_1,...,v_{k_n}).
\end{align}
From the above definitions, we see that $\widehat{c_{t}}^{n} = UU^{\mathsf{T}}$ and $\widetilde{c_{t}}^{n} = VV^{\mathsf{T}}$.
According to Lemma 2.7 in \cite{B1999}, we have
\begin{align}\label{L^4}
	(L(F^{\widehat{c_{t}}^{n}}, F^{\widetilde{c_{t}}^{n} }))^4 &\leq \frac{2}{p^2} \tr ((u_1,...,u_{k_n})(u_1,...,u_{k_n})^{\T})\cdot \tr( \widehat{c_{t}}^{n} + \widetilde{c_{t}}^{n}),
\end{align}
where $L(F,G)$ is the L$\acute{\text{e}}$vy distance between the two probability distribution functions $F$ and $G$. 
For $j=1,...,k_n$ and $i=1,...,p$, under Assumption \ref{asu:dgp} and since $p/k_n = O(1)$, we have
\begin{align*}
	\E[(u_j^{(i)})^2] &=\E\left[ \left( \frac{\int_{(\lfloor tn \rfloor+j-1)/n}^{(\lfloor tn \rfloor+j)/n}b_s^{(i)}ds + \sum_{k=1}^{q}\int_{(\lfloor tn \rfloor+j-1)/n}^{(\lfloor tn \rfloor+j)/n} (\sigma_s^{(i,k)} - \sigma_t^{(i,k)})dB_s^{(k)} }{\sqrt{k_n/n}} \right)^2\right]\\
&\leq \frac{C}{k_n n} + \frac{Cn}{k_n}\cdot  \E\left[\left( \sum_{k=1}^{q}\int_{(\lfloor tn \rfloor+j-1)/n}^{(\lfloor tn \rfloor+j)/n} (\sigma_s^{(i,k)} - \sigma_t^{(i,k)})dB_s^{(k)} \right)^2 \right] \\
(&\text{By Cauchy-Schwarz inequality and It$\hat{\text{o}}$ Isometry}.)\\
& \leq \frac{C}{p n} + \frac{Cp^{2\gamma-1} q}{n^{2\gamma}},
\end{align*}
and similarly
\begin{align*}
	\E[(v_j^{(i)})^2] &= \E \left[ \left( \frac{\sum_{k=1}^{q} \sigma_t^{(i,k)} \Delta_{\lfloor tn \rfloor+j}^nB^{(k)}}{\sqrt{k_n/n}} \right)^2 \right] \leq \frac{Cn}{k_n} \cdot \E\left[ \left( \sum_{k=1}^{q} \Delta_{\lfloor tn \rfloor+j}^nB^{(k)} \right)^2 \right]  \leq \frac{Cq}{p}.
\end{align*}
Based on the above results, we can further obtain
\begin{align}\label{tr:UV}
	\E \left[\tr ((u_1,...,u_{k_n})(u_1,...,u_{k_n})^{\T}) \right] \leq C\left(\frac{p}{n} + pq \left(\frac{p}{n}\right)^{2\gamma}\right), \quad 
	\E\left[\tr( \widehat{c_{t}}^{n} + \widetilde{c_{t}}^{n})\right] = C(pq). 
\end{align}
Plugging them into \eqref{L^4} results in $\E \left[L(F^{\widehat{c_{t}}^{n}}, F^{\widetilde{c_{t}}^{n} })\right] \leq C(\frac{q}{n} + q^2\left(\frac{p}{n}\right)^{2\gamma} ) \rightarrow 0$, if condition $(\romannumeral1)$ is satisfied. Thus, Markov's inequality implies that $L(F^{\widehat{c_{t}}^{n}}, F^{\widetilde{c_{t}}^{n} })$ converges to 0 in probability, namely, the empirical spectral distributions of $\widehat{c_{t}}^{n}$ and $\widetilde{c_{t}}^{n}$ converge to the same distribution $F^{\bar{p}, H_t}$.

On the same probability space $(\Omega, \mathcal{F}, (\mathcal{F}_t)_{0\leq t\leq 1}, \mathbf{P})$ where $X$ is defined, we define, for $p=1,2,...$ and $1 \leq l\leq k_n$, $\mathbf{Z}_{l}^{(p)} = (Z_l^{(p,j)})_{1\leq j\leq p}$ with $Z_l^{(p,j)}$ being $i.i.d.$ standard normal random variables and independent of $B$. Let $c_t^{1/2}$ be the (nonnegative) square root matrix of $c_t$ and
\begin{align}\label{S_p}
S_p := 1/k_n \times \sum_{l=1}^{k_n} c_t^{1/2} \mathbf{Z}_l^{(p)} (\mathbf{Z}_l^{(p)})^\T c_t^{1/2}. 
\end{align}
Condition on $\mathcal{F}_t$, under which $\sigma_t, c_t$ are measurable, we have $\widetilde{c_{t}}^{n}  \overset{d}= S_p$ since their components are $i.i.d.$ and satisfy $\sigma_{t} \sqrt{n}\Delta_i^n B \overset{d}= c_t^{1/2} \mathbf{Z}_l^{(p)} \sim N(\mathbf{0}_{p},c_t)$, thus the ESD $F^{S_p}$ converges to the same limiting spectral distribution of $\widetilde{c_{t}}^{n}$, that is $F^{\bar{p},H_t}$.
According to Proposition \ref{pro1}, the spectral distribution function $F^{\bar{p},H_t}$ is driven by \eqref{equ:st} and \eqref{equ:esd}, with $H$ and $F^{y,H}$ replaced by $H_t$ and $F^{\bar{p},H_t}$ respectively. The proof of the theorem is finished.  $\hfill \square$

\begin{rmk}\label{rmk:lev}
In \cite{ZL2011}, they further require that, there exists a sequence $\eta(p) = o(p)$ and a sequence of index sets $\mathcal{I}_p$ satisfying $\mathcal{I}_p \subset{1,...,p}$ and $\# \mathcal{I}_p \leq \eta_p$ such that the volatility process $\sigma$ may be dependent on $B$, but only on $\{B^{(j)}: j\in \mathcal{I}_p \}$. This condition indeed excludes the general leverage effect in high dimension, and is not really relevant in financial econometrics context. We do not have such a limitation. This is because, with the continuity condition \eqref{asu:dgp:cont} and the required asymptotic conditions, the population covariance matrix of $\Delta_i^n X$ with $i= \lfloor tn \rfloor +1,...,\lfloor tn \rfloor +k_n$ in \eqref{spotvolest} can be approximated by a common covariance matrix $c_{t}$, which is assumed to be c$\grave{\text{a}}$dl$\grave{\text{a}}$g. As a result, $c_{t}$ can be seen as a nonrandom matrix at time $t$ thus is independent of $\Delta_i^nB$ for $i\geq  \lfloor tn \rfloor +1$. 
\end{rmk}

To see how fast can the dimensionality $p$ increases, we now take a deep analysis on condition $(\romannumeral1)$.
Without loss of generality, we assume $p= O(n^{a}), q= O(n^{b})$ with $0< a\leq 1, 0\leq b <1$, then the condition turns to requiring $b+a\gamma - \gamma < 0$. We conclude that: 
\begin{enumerate}
\item $\gamma \neq 0$, otherwise, the condition  can never be satisfied. 
\item  For fixed $\gamma$, we require $\gamma>b$ and $a < 1-b/\gamma$, namely, the fastest divergent rate $p$ can achieve decreases linearly with the slop $1/\gamma$ as $b$ increases. 
\item If $b=0$, corresponding to constant factor number $q$, we only require $a<1$. As a result, almost all of the increments $\Delta_i^n X$ can be used for the spot volatility matrix estimator \eqref{spotvolest}, since $k_n = O(p)$.
\item We allow $q$ increases with $p$ by setting $b>0$. For example, if $p=q$, as considered in \cite{ZL2011}, the condition is $b=a<\frac{\gamma}{\gamma + 1}$. As $\gamma$ increases to $\infty$, the divergence rate allowed for $p$ and $q$ increases and approaches to $O(n)$. Conversely, as $\gamma$ decreases to $0$, the upper bound of $a$ and $b$ also decreases to 0. 
\end{enumerate}

\subsection{Linear spectral statistics (LSS)}\label{sec:lss}
For any random symmetric and nonnegative definite matrix $M$, its linear spectral statistics (LSS) is defined as 
\begin{align}
\frac{1}{p} \sum_{i=1}^{p} f(\lambda_i^{M}) =  \int f(x)dF^{M}(x),
\end{align}
where we recall that $\{\lambda_i^{M}: i=1,..., p\} $ is the real spectrum of $M$, $F^{M}$ is the ESD of $M$, and $f(\cdot)$ is a function defined on $[0,\infty)$.
LSS is important in multivariate statistical inference since many statistics for the population parameters can be expressed in such a form (\cite{A2003}). 
With the definitions and notation given in Proposition \ref{pro1}, a general LSS based on the sample covariance matrix $S_p$ can then be written as $\hat{\theta} = \int f(x)dF^{S_p}(x)$. Proposition \ref{pro1} describes the asymptotic distributional behavior of the eigenvalues of the sample covariance matrix, thus the point-wise limits of the eigenvalue statistics are the integrals of the corresponding functions with respect to the limiting spectral distribution $F^{y, H}$ in \eqref{equ:esd}. Namely, we have $\hat{\theta} \rightarrow \theta:= \int f(x)dF^{y, H}(x)$, almost surely, which serves as a first-order convergence result for LSS.
Furthermore, for making statistical inferences or conducting hypothesis testing problems on the population parameter, we need the second-order convergence result of CLT for LSS. 
To this end, we let $F^{y_n, H_p}$ be the distribution defined by \eqref{equ:st} and \eqref{equ:esd} with the parameters $y, H$ replaced by $y_n:=p/n, H_p:= F^{\Sigma_p}$, respectively, and define $G_p:= p( F^{S_p} - F^{y_n, H_p} )$. Our target is investigate the fluctuation of the following scaled and centralized LSS:
\begin{align}\label{equ:lss}
\int f(x)dG_p(x).
\end{align}
We remark that $G_p$ is defined as the difference between the sample-based ESD $F^{S_p}$ and $F^{y_n, H_p}$, instead of $F^{y, H}$. This is because, as explained in \cite{BS2004}, on one hand, the convergence rate of $y_n\rightarrow y$ and $H_p \rightarrow H$ can be arbitrarily slow; on the other hand, from the perspective of statistical inference, $H_p$ can be viewed as a description of the current population and $y_n$ is the ratio of dimension to sample size for the current sample, such a consideration is more realistic. The empirical result in the large dimensional random matrix field shows that \eqref{equ:lss} converges to a Gaussian distribution, whose expectation and covariance are obtained by using the Stieltjes transform. 
\begin{pro}[Theorem 1.1 in \cite{BS2004}]\label{pro:lss}
Assume on a common probability space, we have
\begin{enumerate}[label=(\roman*)]
	\item for $p=1,2,...$ and for $1 \leq l\leq n$, $\mathbf{Z}_{l}^{(p)} = (Z_l^{(p,j)})_{1\leq j\leq p}$ with $Z_l^{(p,j)}$ $i.i.d.$ with mean 0, variance 1, and finite fourth moment;
	\item as $n\rightarrow \infty$, it holds that $p\rightarrow \infty$ with $y_n:= p/n \rightarrow y>0$;
	\item $\Sigma_p$ is a nonrandom Hermitian nonnegative definite $p \times p$ matrix such that its spectral norm is bounded in $p$ and its ESD $H_p:= F^{\Sigma_p}$ converges almost surely to a probability distribution $H$ on $[0, \infty)$ as $p \rightarrow \infty$.
\end{enumerate}
Let $\Sigma_p^{1/2}$ be the (nonnegative) square root matrix of $\Sigma_p$ and $S_p := 1/n \times \sum_{l=1}^{n} \Sigma_p^{1/2} \mathbf{Z}_l^{(p)} (\mathbf{Z}_l^{(p)})^{\T} \Sigma_p^{1/2}$, $f_1,...,f_k$ are functions defined on $\mathbb{R}$ and are analytical on an open interval containing 
\begin{align*}
\left[ \lim\inf_p \lambda_{\min}^{\Sigma_p} \mathbf{I}_{(0,1)}(y) (1-\sqrt{y})^2,   \lim\sup_p \lambda_{\max}^{\Sigma_p} (1+\sqrt{y})^2 \right],
\end{align*}
then, with $G_p := p( F^{S_p} - F^{y_n, H_p} )$, the random vector 
\begin{align}\label{lss:rv}
\left(\int f_1(x)G_p(x),..., \int f_k(x)G_p(x)\right)
\end{align}
forms a tight sequence in $p$. Furthermore, if $Z_l^{(p,j)}$ and $\Sigma_p$ are real and $\E[|Z_l^{(p,j)}|^4] = 3$, then \eqref{lss:rv} converges in distribution to a Gaussian vector $(X_{f_1},...,X_{f_k})$ with means
\begin{align}\label{lss:mean}
\E [X_f] = - \frac{1}{2\pi i} \int f(z) \frac{y\int \underline{m}(z)^3 t^2 (1+t\underline{m}(z))^{-3} dH(t)}{(1-y\int \underline{m}(z)^2 t^2 (1+t\underline{m}(z))^{-2} dH(t))^2} dz
\end{align}
and covariance function 
\begin{align}\label{lss:cov}
Cov(X_f, X_g) = -\frac{1}{2\pi^2} \int \int \frac{f(z_1)g(z_2)}{(\underline{m}(z_1)-\underline{m}(z_2))^2} \frac{d}{dz_1}\underline{m}(z_1) \frac{d}{dz_2} \underline{m}(z_2)dz_1dz_2
\end{align}
where, $f,g \in \{f_1,...,f_k\}$ and we recall that  $\underline{m}(z)$ is the Stieltjes transform of $\underline{F}^{y,H} = y F^{y,H} + (1-y) \mathbf{1}_{[0, \infty)}$, which is the limiting empirical distribution function of $\underline{S}_p := 1/n \times \sum_{l=1}^{n} (\mathbf{Z}_l^{(p)})^{\T} \Sigma_p \mathbf{Z}_l^{(p)}$. The contours in \eqref{lss:mean} and \eqref{lss:cov} (two in \eqref{lss:cov}, which we may assume to be nonoverlapping) are closed and are taken in the counter-clockwise direction in the complex plane, each enclosing the support of the limiting spectral distribution $F^{y, H}$.
\end{pro}

As a direct application of Proposition \ref{pro:lss}, under the special case of constant volatility model \eqref{constvol}, after directly replacing $\Sigma_p$ and $S_p$ with the integrated covariance matrix $\int_{0}^{1} c_td_t$ and its realized estimator $\widehat{IV}^n$ in \eqref{intvolest} respectively, we see that the random vector defined as in \eqref{lss:rv} will converge to a Gaussian distribution with expectation and covariance given in \eqref{lss:mean} and \eqref{lss:cov}, if $y_n:= p/n \rightarrow y>0$ as $n\rightarrow \infty$. 

Before presenting our LSS result for the spot volatility matrix, we make some remarks as follows. 
We note that, in the above proposition, the population covariance matrix $\Sigma_p$ is restricted to be nonrandom, while such a limitation is not required for the first-order limiting result in Proposition \ref{pro1}. 
We do not need such a restrictive nonrandom condition for the volatility process $\sigma$ in \eqref{model}, even for the investigation of LSS for the spot covariance matrix. 
Similar to the analysis in Remark \ref{rmk:lev}, we can approximate the population covariance matrix of $\Delta_i^n X$ with $i= \lfloor tn \rfloor +1,...,\lfloor tn \rfloor +k_n$ in \eqref{spotvolest} by a common covariance matrix $c_{t}$, which can be seen as a nonrandom matrix at time $t$. 
Thus, the nonrandom assumption for $\sigma$ is not necessary for us.
 
For the spot volatility process in model \eqref{model}, we define, for fixed $t \in[0,1]$, 
\begin{align}
\widehat{G}_{p,t} = p (F^{\widehat{c_{t}}^{n}} - F^{z_n, H_{p,t}}),
\end{align}
where $z_n:=p/k_n$ and $H_{p,t} := F^{c_t}$, whose definition is given in Theorem \ref{spot_esd}.
We have

\begin{thm}\label{thm:lss}Under Assumption \ref{asu:dgp} and the following conditions:
\begin{enumerate}[label=(\roman*)]
	\item as $n \rightarrow \infty$, $p \rightarrow \infty$ and meanwhile  $q/n \rightarrow 0$,  $p^{3/2}/n \rightarrow 0$, $qp^{2\gamma+3/2}/n^{2\gamma} \rightarrow 0$ and $z_n:=p/k_n \rightarrow \bar{p}>0$;
	\item For the $p \times p$ random, symmetric and nonnegative definite matrix $c_t$ with $t\in[0,1]$, its spectral norm is bounded in $p$ and its ESD $H_{p,t}:= F^{c_t}$ converges almost surely to a probability distribution $H_t$ on $[0, \infty)$ as $p \rightarrow \infty$;
\end{enumerate}
Let $f_1,...,f_k$ be functions on $\mathbb{R}$ and are analytical on an open interval containing 
\begin{align*}
	\left[ \lim\inf_p \lambda_{\min}^{c_t} \mathbf{I}_{(0,1)}(\bar{p}) (1-\sqrt{\bar{p}})^2,   \lim\sup_p \lambda_{\max}^{c_t} (1+\sqrt{\bar{p}})^2 \right],
\end{align*}
then, condition on $\mathcal{F}_t$ and with $\widehat{G}_{p,t} =p ( F^{\widehat{c_{t}}^{n}} - F^{z_n, H_{p,t}} )$, the random vector 
\begin{align}\label{lss:spot}
	\left(\int f_1(x)\widehat{G}_{p,t}(x),..., \int f_k(x)\widehat{G}_{p,t}(x)\right)
\end{align}
forms a tight sequence in $p$ and 
converges in distribution to a Gaussian vector $(X_{f_1},...,X_{f_k})$ with means and covariance function given in \eqref{lss:mean} and \eqref{lss:cov} respectively.
\end{thm}
\textit{Proof.} Recall that $\widehat{c_{t}}$ and $\widetilde{c_{t}}$ are defined as in \eqref{spotvolest} and \eqref{def:sigtilde}, we further define $\widetilde{G}_{p,t} = p( F^{\widetilde{c_{t}}^{n}} - F^{z_n, H_{p,t}} )$ and will show that
\begin{align}\label{thm2:op1}
	\int f(x)d\widehat{G}_{p,t} - 	\int f(x)d\widetilde{G}_{p,t} = o_p(1).
\end{align}
With $U$ and $V$ defined as in \eqref{def:UV}, we can obtain
\begin{align}\label{thm2:s0}
	\begin{split}
		&\left | \int f(x)d\widehat{G}_p - 	\int f(x)d\widetilde{G}_p  \right | = \left | \sum_{i=1}^{p} f(\lambda_i^{\widehat{c_{t}}^{n}}) -  f(\lambda_i^{\widetilde{c_{t}}^{n}}) \right| \\
		& \leq C \sum_{i=1}^{p} \left| \lambda_i^{\widehat{c_{t}}^{n}} - \lambda_i^{\widetilde{c_{t}}^{n}} \right| \\
		& \text{(According to H$\ddot{\text{o}}$lder's inequality)} \\
		& \leq C \sqrt{\sum_{i=1}^{p} \left( \sqrt{\lambda_i^{\widehat{c_{t}}^{n}}} + \sqrt{\lambda_i^{\widetilde{c_{t}}^{n}}} \right)^2 } \sqrt{\sum_{i=1}^{p} \left( \sqrt{\lambda_i^{\widehat{c_{t}}^{n}}} - \sqrt{\lambda_i^{\widetilde{c_{t}}^{n}}} \right)^2 }  \\
		&(\text{According to Theorem A.37 in \cite{B2010}})\\
		&\leq C \sqrt{2p\cdot (\lambda_{\max}^{\widehat{c_{t}}^{n}} + \lambda_{\max}^{\widetilde{c_{t}}^{n}})} \tr( (U-V)(U-V)^{\T}).
	\end{split}
\end{align}
By using the results in \eqref{tr:UV}, we have 
\begin{align}\label{thm2:s1}
	\|U-V\|_F^2 =  \tr (U-V)(U-V)^{\T} = O_p\left(\frac{p}{n} + pq \left(\frac{p}{n}\right)^{2\gamma}\right).
\end{align}
According to \cite{YBK1988}, we have, almost surely,
\begin{align}\label{thm2:s2}
	\lambda_{\max}^{\widetilde{c_{t}}^{n}}  \leq \lim\sup_{p} \|c_t\| (1+\sqrt{\bar{p}})^2<\infty,
\end{align}
and since $\widehat{c_{t}}^{n} = UU^{\mathsf{T}}$ and $\widetilde{c_{t}}^{n} = VV^{\mathsf{T}}$, we have 
\begin{align}\label{thm2:s3}
	\lambda_{\max}^{\widehat{c_{t}}^{n}} = \|U\|_2^2 \leq (\|V\|_2 + \|U-V\|_2)^2 \leq \lambda_{\max}^{\widetilde{c_{t}}^{n}} + C \|U-V\|_{F} < \infty.
\end{align} 
After plugging \eqref{thm2:s1}--\eqref{thm2:s3} into \eqref{thm2:s0} and using condition $(\romannumeral1)$, we can obtain \eqref{thm2:op1}. Similarly, by defining $S_p$ as in \eqref{S_p} and following the analyses therein, directly applying Proposition \ref{pro:lss} yields that the random vector $\left(\int f_1(x)\widetilde{G}_{p,t}(x),..., \int f_k(x)\widetilde{G}_{p,t}(x)\right)$
forms a tight sequence in $p$ and 
converges in distribution to a Gaussian vector $(X_{f_1},...,X_{f_k})$, whose means and covariance function are given by \eqref{lss:mean} and \eqref{lss:cov} respectively. And \eqref{thm2:op1} implies that the same result holds for the random vector \eqref{lss:spot}, which finishes the proof of the theorem.  $\hfill \square$

To the best of our knowledge, allowing every entry of $\Sigma_{p}$ in Proposition \ref{pro:lss} to be random and extending related theory has not been considered, even in the large dimensional random matrix community. Partial of this problem have been solved if $\Sigma_{p}$ has some extra structures. For example, for the elliptical correlated model with $\Sigma_{p} = w \Lambda_p$, where $w>0$ is a scalar random variable and $\Lambda_p$ is a $p\times p$ nonrandom matrix of full rank (See, e.g. \cite{HLLZ2019}, \cite{HLLZ2025} and references therein.).
In \cite{YZC2021}, they considered the elliptically distributed samples, which are similar to the high-frequency increments under the class $\mathcal{C}$ of \eqref{vol:classC}. They established the central limit theorem for the LSS of the sample covariance matrix by using the self-normalized observations, such an estimator was proposed in \cite{ZL2011} and its limiting spectral distribution was also investigated therein.
In Theorem \ref{thm:lss}, we do require the mentioned structure for the volatility process $c$ and all of its entries can be random. This also inspires us that, establishing the LSS result similar to Proposition \ref{pro:lss} for the sample covariance matrix when the entries of $\Sigma_{p}$ are random is also possible, at least it can be realized by restricting the variance of the entries. 

Similar to the analyses for Theorem \ref{spot_esd}, we assume $p= O(n^{a})$ with $0< a < 2/3$ and $q= O(n^{b})$ with $0\leq b <1$, then condition $(\romannumeral1)$ in Theorem \ref{thm:lss} turns to requiring $b+3a/2+2a\gamma - 2\gamma < 0$. We conclude that:
\begin{enumerate}
	\item $\gamma \neq 0$, otherwise, the condition  can never be satisfied. 
	\item  For fixed $\gamma$, we require $2\gamma>b$ and $a < (2\gamma-b)/(3/2+2\gamma)$, namely, the fastest divergent rate $p$ can achieve decreases linearly with the slop $1/(2\gamma)$ as $b$ increases. 
	\item If $b=0$, corresponding to constant factor number $q$, we require $a<2\gamma/(3/2+2\gamma)$. When $\gamma>3/2$, $p$ can approach the possible fastest divergence rate of $n^{2/3}$, while the divergence rate decreases to 0 as $\gamma$ tends to 0. 
	\item It is allowed that $q$ increases with $p$ by setting $b>0$. For example, if $p=q$, the condition is $a=b<\frac{2\gamma}{2\gamma + 5/2}$. As $\gamma$ increases to $\infty$, the divergence rate allowed for $p$ and $q$ increases and approaches to $O(n^{2/3})$. Conversely, as $\gamma$ decreases to $0$, the upper bound of $a$ and $b$ also decreases to 0. 
\end{enumerate}

\begin{rmk}
\red{Both in Theorem \ref{spot_esd} and Theorem \ref{thm:lss}, the sample size $k_n$ used for the spot volatility matrix estimator \eqref{spotvolest} depends on the divergence rate of $p$ via $p/k_n\rightarrow \bar{p}>0$, based on which $k_n$ should be selected. From the analyses of the asymptotic conditions in these theorems, we see that the divergence rate allowed for $p$ is affected by $\gamma$ in \eqref{asu:dgp:cont} and the factor number $q$. In practice, the factor number is often treated to be constant, as shown in \cite{FF1993}, \cite{FF2015}, and some others.
As for $\gamma$, it varies for different assets and is also time-varying, since it relies on the variational pattern of the volatility process. We may get a rough knowledge about $\gamma$ via estimating the spot volatility matrix process by using the historical data and investigating the continuity of the estimated path. 
With these considerations, the selection of $k_n$ can be determined by satisfying the asymptotic condition $(\romannumeral1)$ in Theorem \ref{spot_esd} and Theorem \ref{thm:lss}.}
\end{rmk}


\section{Application}\label{sec:app}
In this section, we apply the theoretical results established in the last section to conduct hypothesis testing problems of identity test and sphericity test for the spot volatility matrix $c_t$ for $t\in[0,1]$.

\subsection{Identity test}\label{sec:indtest}
We first consider testing if the spot volatility matrix $c_{t}$ is equal to a given matrix $\Sigma$ or not. Assume $\Sigma$ is invertible, we note that such a testing problem is equivalent to the identity test after multiplying the observations by $\Sigma^{-1/2}$. Thus, without loss of generality, we consider the following hypothesis testing problem for some $t\in[0,1]$:
\begin{align}\label{testing}
	H_0: c_{t} = \mathbb{I}_{p}  \qquad \text{v.s.} \qquad  H_1: c_{t} \neq \mathbb{I}_{p}.
\end{align}

Before giving out our test statistic, we review the related results for the empirical covariance matrix under fixed dimension first. Suppose that $\bm{x}$ follows a $p$-dimensional multivariate Gaussian distribution $N(\mu_p, \Sigma_p)$ and we want to test if $\Sigma_p  = \mathbb{I}_p$ or not. Let $(\bm{x}_1,...,\bm{x}_n)$ be a set of mutually independent samples from $\bm{x}$, then the likelihood ratio test (LRT) statistic is defined as 
\begin{align}\label{LRTstat}
\tr(S_n)  - \log(|S_n|) - p,
\end{align}
where 
\begin{align}\label{S_n}
S_n:= \frac{1}{n} \sum^{n}_{i=1}(\bm{x}_i - \bar{\bm{x}})(\bm{x}_i - \bar{\bm{x}})^{\T}
\end{align}
 is the sample covariance matrix with $\bar{\bm{x}} = \frac{1}{n} \sum_{i=1}^{n} \bm{x}_i$. For fixed $p$, as $n \rightarrow \infty$, the classical theory from the multivariate statistical analysis states that \eqref{LRTstat} converges to $\chi_{1/2p(p+1)}^2$ distribution if $\Sigma_p  = \mathbb{I}_p$ (\cite{A2003}). But under the large dimensional regime, such a test statistic is infeasible, as justified by the empirical studies in \cite{BJYZ2009}. Moreover, based on the central limit theorem results for LSS in Proposition \ref{pro:lss}, \cite{BJYZ2009} studied the statistic \eqref{LRTstat} and demonstrated that \eqref{LRTstat} converges to a normal distribution when $p/n\rightarrow y \in (0,1)$ as $n\rightarrow \infty$\footnote{\red{For the critical case of $y=1$, \cite{JJY2012} and  \cite{JY2013} studied the statistic \eqref{LRTstat} by using the Selberg integral and analyzing its moments respectively, their results demonstrate that \eqref{LRTstat} still converges to a normal distribution, but the expectation and variance are different with the ones obtained for $y\in (0,1)$. 
If $y\in (1,\infty)$, the statistic \eqref{LRTstat} is degenerate and not applicable since $|S_n| = 0$ when $p>n$, thus the statistics has to be modified. To this end, \cite{QLY2023} proposed a quasi-LRT test statistic and investigated its limiting distribution based on the LSS of a special estimator called re-normalized sample covariance matrix.}}. 
 Similarly, applying our theoretical results in Theorem \ref{thm:lss}, we can propose a test statistic for testing the hypotheses in \eqref{testing}, as follows. 
\begin{thm}\label{thm:testing}
Define $g(x) = x-\log(x) - 1$, $ v(g) = -2\log (1-\bar{p}) - 2\bar{p}$ and $m(g) = -\frac{\log(1-\bar{p})}{2}$, and 
\begin{align}
\widehat{L}_n = \tr(\widehat{c_{t}}^{n} ) - \log(|\widehat{c_{t}}^{n}|) - p.
\end{align} 
Under the conditions of Theorem \ref{thm:lss}, for $\bar{p} \in (0,1)$ and under $H_0$, we have 
\begin{align}\label{thm:test:res}
(v(g))^{-1/2} (\widehat{L}_n- p\cdot F^{z_n}(g) - m(g))\longrightarrow^d N(0,1),
\end{align}
where $F^{z_n}$ is the Mar$\check{\text{c}}$enko-Pastur law with ratio index $z_n$ and scale index $1$, namely $F^{z_n}(g) = 1+(\frac{1}{z_n}-1)\log(1-z_n)$.
\end{thm}
\textit{Proof.} According to the proof of Theorem \ref{thm:lss}, almost surely, $\widehat{c_{t}}^{n}$ has the same asymptotic distribution with $S_p$ defined in \eqref{S_p}. We define 
\begin{align*}
	\widetilde{L}_n = \tr(S_p ) - \log(|S_p|) - p. 
\end{align*}
Under the null hypothesis $H_0$ of $c_t = \mathbb{I}_p$, we have
\begin{align*}
	\widetilde{L}_n &= \tr(S_p ) - \log(|S_p|) - p\\
	& = \sum_{i=1}^{p} (\lambda_i^{S_p} - \log (\lambda_i^{S_p}) - 1) \\
	&= p\cdot \int g(x) d(F^{S_p}- F^{z_n, H_{p,t}} )) + p\cdot F^{z_n}(g), 
\end{align*}
and the first term converges to a Gaussian vector with mean $m(g)$ and variance $v(g)$, which are calculated in the proof of Theorem 3.1 in \cite{BJYZ2009}. A direct application of Slutsky's theorem yields the conclusion \eqref{thm:test:res}. $\hfill \square$

Conclusion in Theorem \ref{thm:testing} gives us the following test statistic:
\begin{align}
\text{BJYZ-test-spot}: \quad \frac{ \widehat{L}_n- p\cdot \left(1+(\frac{k_n}{p}-1)\log(1-p/k_n)\right)- m(g)} {(v(g))^{1/2}}.
\end{align}

Noticing that in the above theorem (Also for the results in \cite{BJYZ2009}), the limiting value of the ratio between the dimension $p$ and the sample size $k_n$, $\bar{p}$, is constrained to be smaller than 1, which may not be satisfied for some applications in practice. 
When the dimension $p$ can be larger than the sample size $n$ in the large dimensional regime, under the $i.i.d.$ setting, \cite{LW2002} proposed the following test statistic:
\begin{align*}
\frac{1}{p} \tr((S_n - \mathbb{I}_p)^2) - \frac{p}{n} \left(\frac{1}{p} \tr(S_n)\right)^2+\frac{p}{n},
\end{align*}
and established its asymptotic properties. 
Similarly, for the spot volatility matrix estimator using the high-frequency data, we can define 
\begin{align}\label{LW1}
	\widehat{W}_n = \frac{1}{p} \tr((\widehat{c_{t}}^{n} - \mathbb{I}_p)^2) - \frac{p}{k_n} \left(\frac{1}{p} \tr(\widehat{c_{t}})\right)^2+\frac{p}{k_n},
\end{align}
and we can obtain
\begin{thm}\label{thm:testing:lw}
	Under the conditions of Theorem \ref{thm:lss}, with $\bar{p}$ being any real number within $ (0,\infty)$, we have, under $H_0$,  
	\begin{align}\label{thm:test:res2}
	k_n\widehat{W}_n - p \longrightarrow^{d} N(1,4).
	\end{align}
\end{thm}
\textit{Proof.} We define
\begin{align*}
	\widetilde{W}_n = \frac{1}{p} \tr((S_p - \mathbb{I}_p)^2) - \frac{p}{k_n} \left(\frac{1}{p} \tr(S_p)\right)^2+\frac{p}{k_n},
\end{align*}
with $S_p$ defined in \eqref{S_p}. According to \eqref{thm2:op1} in the proof of Theorem \ref{thm:lss}, we have $p(\widetilde{W}_n - \widehat{W}_n) = o_p(1)$, thus $p(\widetilde{W}_n - \widehat{W}_n) \longrightarrow^p 0$. Moreover, it holds that $k_n\widehat{W}_n - p \longrightarrow^{d} N(1,4)$, as given by Proposition 7 in \cite{LW2002}, where the required conditions are satisfied under our setting. This finishes the proof of Theorem \ref{thm:testing:lw}. 
$\hfill \square$

With Theorem \ref{thm:testing:lw}, we can obtain another test statistic 
\begin{align}
 \text{LW-test-spot}: \quad \frac{k_n\widehat{W}_n - p-1}{2}. 
\end{align}

\begin{rmk}
In \cite{WYMC2013}, they noticed that the statistic in \cite{BJYZ2009} only works for the case when the population mean is known, which is quite restrictive in practice and is significantly different from the theory for the unknown population mean case (Refer to \cite{B2010} and \cite{P2014} for more discussion). They solved the mentioned problem by estimating the population mean vector via its sample mean estimator and established the asymptotic properties for the sample covariance matrix estimator. Moreover, the entries of the sample data are required to be of normal distribution in \cite{BJYZ2009} and \cite{LW2002}, \cite{WYMC2013} further relaxed this condition by considering a general value of the fourth moment. We note that these mentioned modifications in \cite{WYMC2013} are not necessary for us, because the high-frequency increments are approximately distributed as normal random variables with zero mean. This point can be seen from the proof of Theorem \ref{thm:lss}.
\end{rmk}
%


\subsection{Sphericity test}\label{sec:sphetest}
In high-frequency financial econometrics, testing if the spot volatility matrix $c_t$ equals to a know covariance matrix may be too restrictive and has limited applications. But knowing the covariance matrix up to a constant can be good enough in some applications. For example, the minimum variance portfolio is given by $\Sigma^{-1}\mathbf{1}_{p}/(\mathbf{1}_{p}^{\T} \Sigma \mathbf{1}_{p})$, where $\Sigma$ is the population covariance matrix of all assets used for the investment. The optimal portfolio is therefore invariant to scaling in the covariance matrix. This inspires us to conduct the following sphericity hypothesis testing problem:
\begin{align}\label{testing:sphe}
	H_0: c_{t} = k^2\mathbb{I}_{p}  \qquad \text{v.s.} \qquad  H_1: c_{t} \neq k^2 \mathbb{I}_{p},
\end{align}
where $k^2$ is an unknown scalar parameter.

For the sphericity test, in traditional multivariate analysis based on the sample covariance matrix in $\eqref{S_n}$, when the dimension $p$ is fixed, the classical likelihood ratio test statistic can be found in Section 8.3.1  in \cite{M1982}, and \cite{J1971} proposed a test statistic by using the self-normalized sample covariance matrix and its spectrum. 
In high-dimensional situations, similarly to the discussion in Section \ref{sec:indtest}, the former one is only valid when the dimension is smaller than the sample size, while the latter one does not have such a restriction. Thus, the statistic in \cite{J1972} is more favorable and has been extensively studied in recent years: See, for example, \cite{LW2002}, \cite{WY2013}, \cite{TLL2015} for the linear transform model, which corresponds to taking $\gamma_t$ in \eqref{vol:classC} as a constant; 
\cite{ZPFW2014}, \cite{PV2016} for the elliptical model if replacing $\gamma_t$ in \eqref{vol:classC} by a random variable independent of time $t$. Similarly, we will modify the test statistic in \cite{J1972} under our high-frequency setting in large dimension. We do not consider other statistics and only note that many other ones can be found in \cite{HLLZ2019}, \cite{YZC2021} and references therein. 

We recall that the test statistic in \cite{J1972} based on the sample covariance matrix \eqref{S_n} is given by
\begin{align*}
\frac{1}{p}\tr\left( \left(\frac{S_n}{(1/p) \tr(S_n)} - \mathbb{I}_p\right)^2 \right). 
\end{align*}
For the spot volatility matrix estimator using the high-frequency data, after replacing $S_n$ with $\widehat{c_{t}}^{n}$, we can obtain
\begin{align}\label{teststa_sph}
\widehat{M}_n = \frac{1}{p}\tr\left(\left(\frac{\widehat{c_{t}}^{n}}{(1/p) \tr(\widehat{c_{t}}^{n})} - \mathbb{I}_p\right)^2\right).
\end{align}
For the above estimator, we have
\begin{thm}\label{thm:testing:sph}
	Under the conditions of Theorem \ref{thm:lss}, with $\bar{p}$ being any real number within $ (0,\infty)$, we have, under $H_0$,  
	\begin{align}\label{thm:testres:sph}
		k_n\widehat{M}_n  - p \longrightarrow^{d} N(1,4).
	\end{align}
\end{thm}
\textit{Proof.} The conclusion can be obtained by following the proof of Theorem \ref{thm:testing:lw} and using Proposition 3 in \cite{LW2002}. $\hfill \square$

With Theorem \ref{thm:testing:sph}, we can obtain the following statistic for the sphericity testing problem:
\begin{align}
 \text{J-test-spot:}\quad   \frac{k_n\widehat{M}_n - p-1}{2}. 
\end{align}

\section{Simulation studies}\label{sec:sim}
In this supplementary material, we conduct some simulation studies to evaluate the finite sample performance of the realized spot volatility matrix estimator and verify the theoretical results established in the paper.  

With the underlying efficient log-price process $\{X_t\}_{0 \leq t\leq 1}$ generated from \eqref{model}, we consider two different scenarios for the time-varying volatility process, deterministic or stochastic, with the following forms: 
\begin{align}\label{det:vol}
	\text{Deterministic volatility: } \sigma_t = \sqrt{0.0009 + \text{r}_1 \cdot \sin(2\pi t)} \cdot  \mathbb{I}_p, 
\end{align}
and 
\begin{align}\label{sto:vol}
	\text{Stochastic volatility: } \sigma_t =\left(\sqrt{0.0009}+ \int_{0}^{t} \tilde{b}_sds + \int_{0}^{t}  \rho\tilde{\sigma}_s dB^{(1)}_s + \int_{0}^{t}  \sqrt{1-\rho^2}\tilde{\sigma}_s dW_s\right) \mathbb{I}_p. 
\end{align}
In the above models, $W$ is a one-dimensional standard Brownian motion independent of $B$, $B^{(1)}$ is the first entry of $B$ and $\rho$ controls the correlation strength between $X$ and $\sigma$, $\tilde{b} \in \mathbb{R}$ and $\tilde{\sigma} \in \mathbb{R}$ are the drift process and volatility process of the volatility process $\sigma$ respectively. 
For simplicity, we let $\tilde{\sigma} \equiv \text{r}_2$, which is constant, and $b, \tilde{b} \equiv 0$, since the drift processes play no roles asymptotically. The deterministic volatility and corresponding parameter setting are also considered in \cite{ZL2011} for experimental study. For the stochastic volatility model, It$\hat{\text{o}}$ semimartingale like $X$ is used and is popular in high-frequency literature. 
We consider different setting for $\text{r}_1 = 0, 0.0004, 0.0008$ in \eqref{det:vol} and $\text{r}_2= 0, 0.01, 0.02$ in \eqref{sto:vol}, these parameters control the variational intensity of the volatility process. 
Without loss of generality, we focus our analysis at the time point $t=0$ and fix $n=4680$, corresponding to 5-second observational data within a 6.5-hour trading day in the real financial market.  We consider $k_n = \lfloor \sqrt{n} \rfloor= 68$ and $p/k_n \equiv \overline{p}$ for different values of $\bar{p} = 0.5 ,1, 1.5$ (Correspondingly, $p=34,68,102$). 

\begin{figure}[!htbp]
	\centering 
	\subfigure[$p=34, k_n=68$]{\includegraphics[width=5cm,height=4cm]{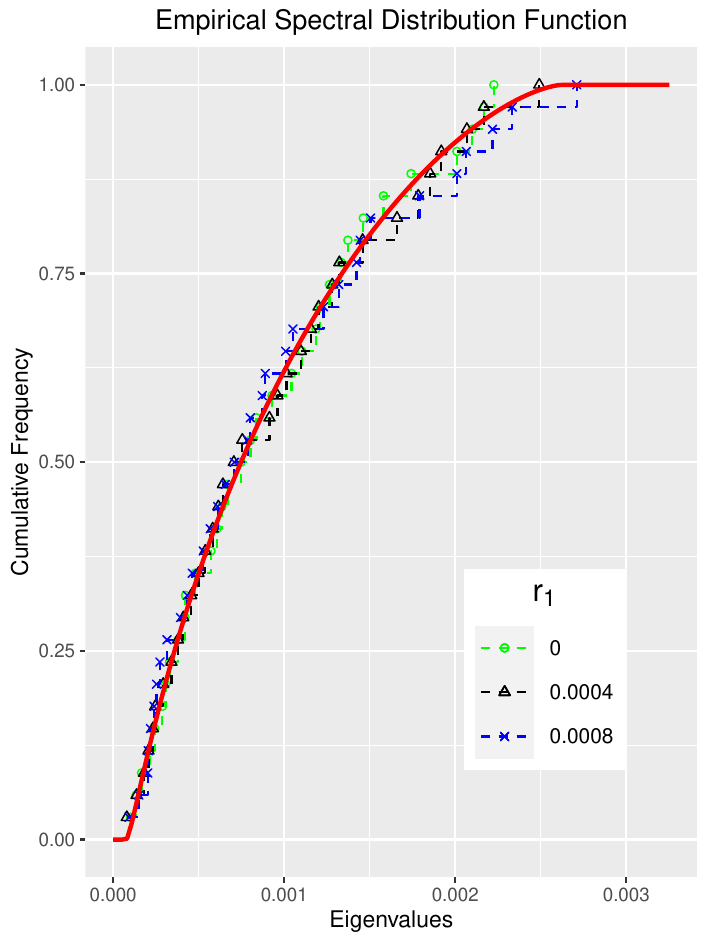}}	
	\subfigure[$p=68, k_n=68$]{\includegraphics[width=5cm,height=4cm]{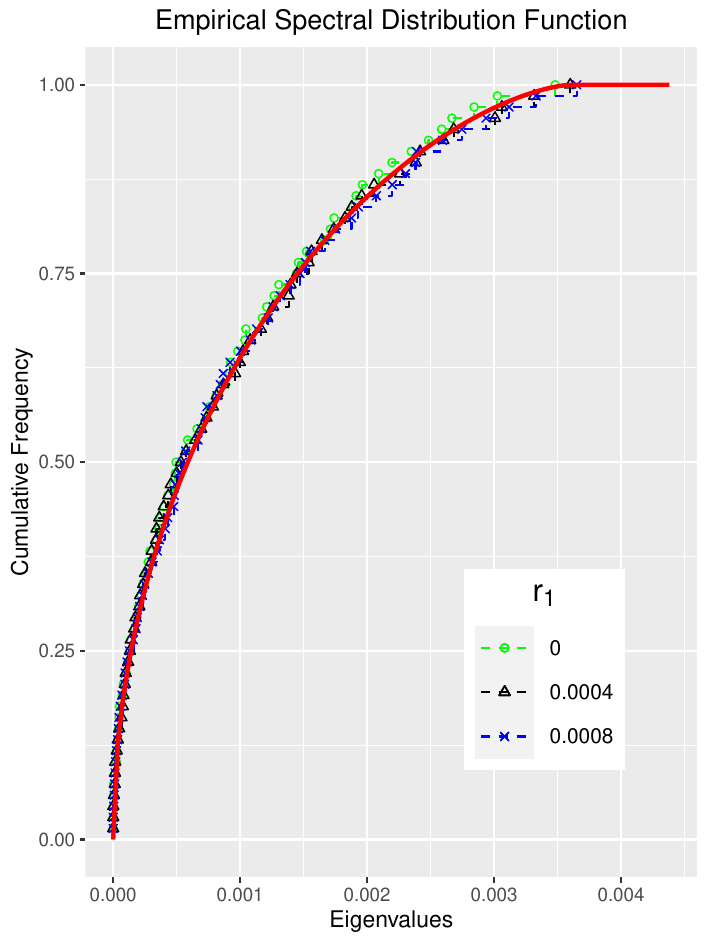}}
	\subfigure[$p=102, k_n=68$]{\includegraphics[width=5cm,height=4cm]{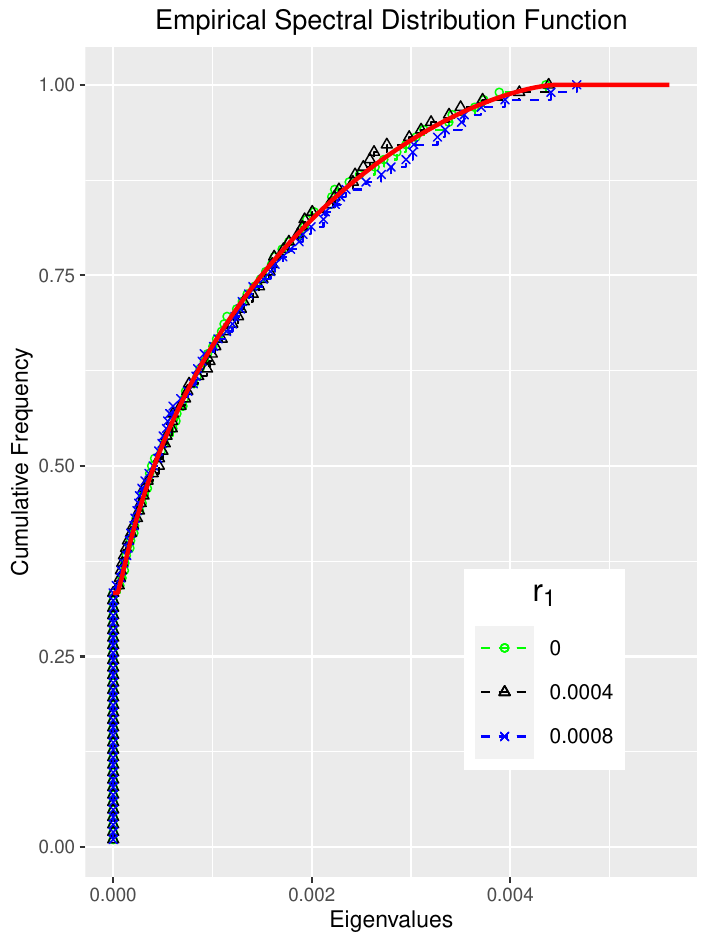}}
	\subfigure[$p=34, k_n=68$]{\includegraphics[width=5cm,height=4cm]{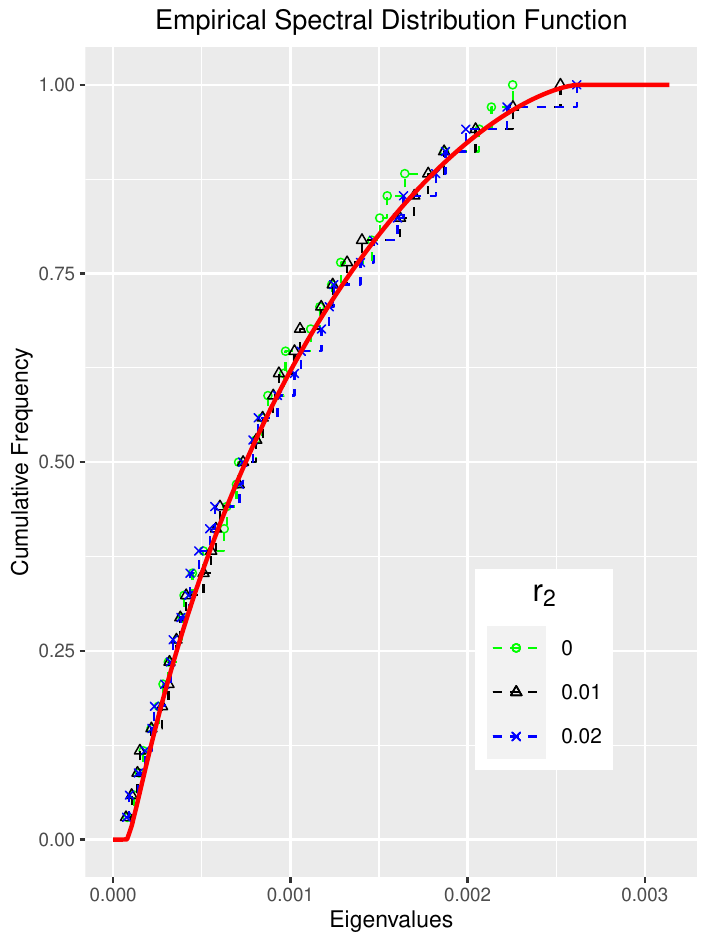}}
	\subfigure[$p=68, k_n=68$]{\includegraphics[width=5cm,height=4cm]{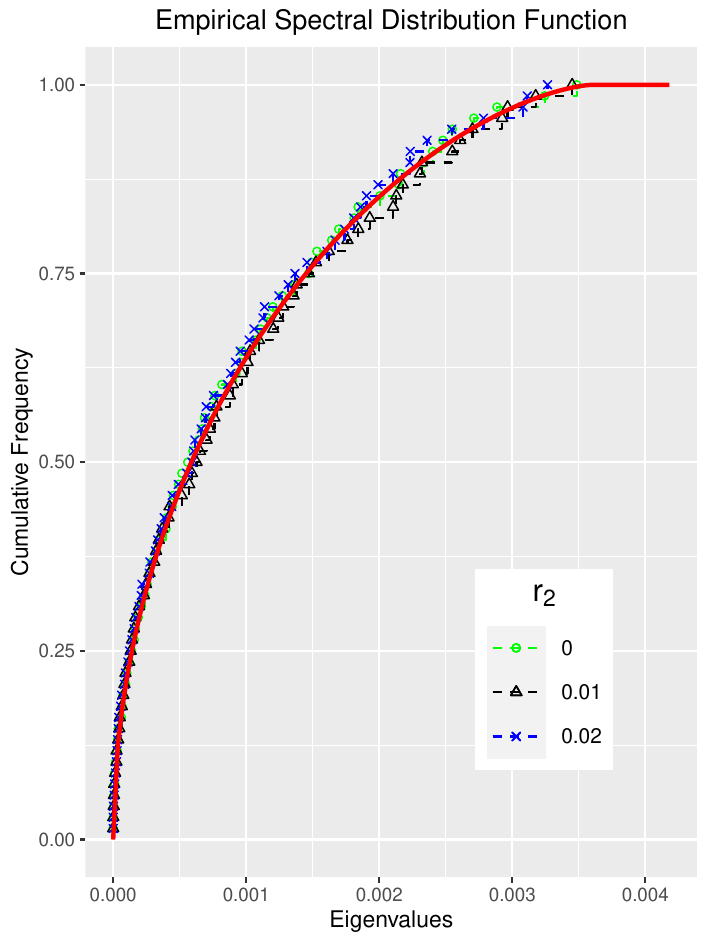}}
	\subfigure[$p=102, k_n=68$]{\includegraphics[width=5cm,height=4cm]{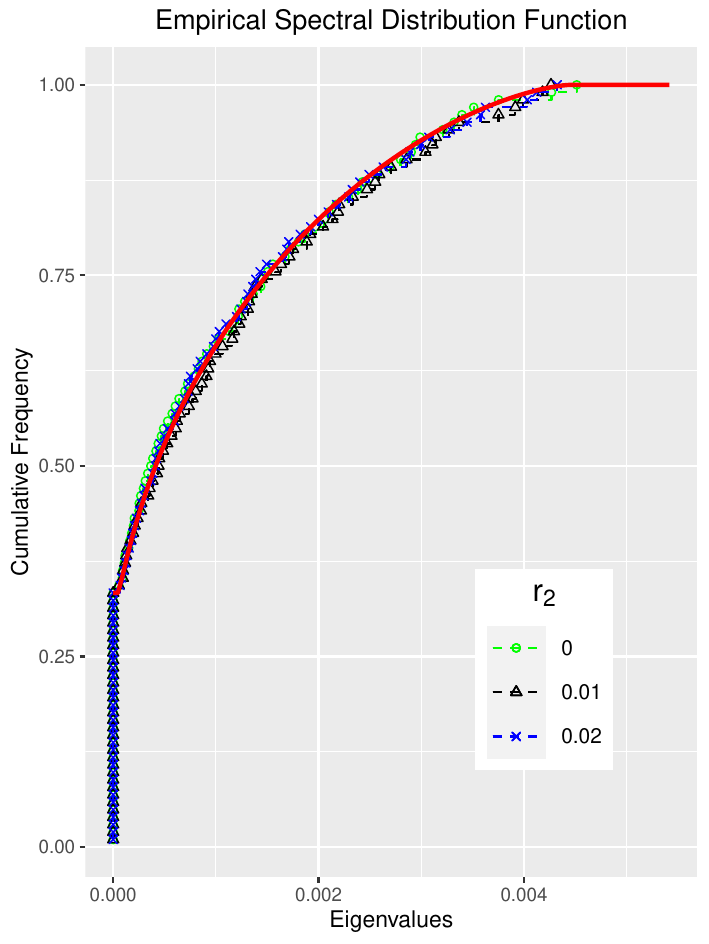}}
	\caption{The empirical spectral distribution functions of realized spot volatility matrix at $t=0$ for the deterministic volatility model \eqref{det:vol} (upper) and the stochastic volatility model \eqref{sto:vol} with $\rho=0$ (lower), with fixed $k_n=68$ and $p=34,68,102$ (Correspondingly, $\bar{p} = 0.5 ,1, 1.5$.). The red solid lines are corresponding theoretical limiting spectral distribution functions, namely the Mar$\check{\text{c}}$enko-Pastur law given by Theorem \ref{spot_esd}.}
	\label{fig_esd}
\end{figure}

\begin{figure}[!htbp]
	\centering 
	\subfigure[$p=34, k_n=68$]{\includegraphics[width=5cm,height=4cm]{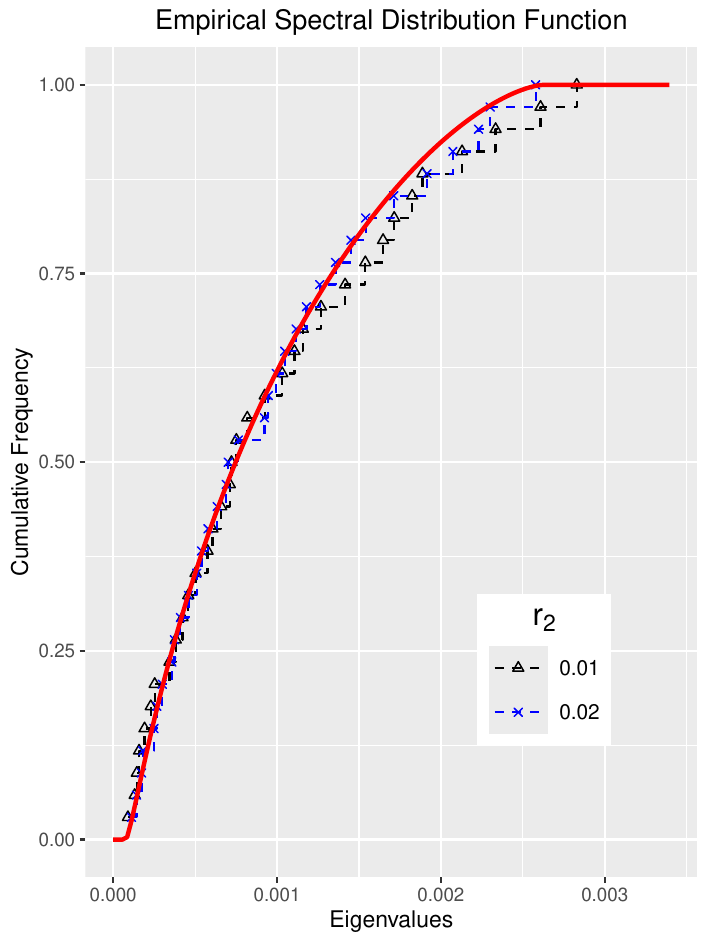}}	
	\subfigure[$p=68, k_n=68$]{\includegraphics[width=5cm,height=4cm]{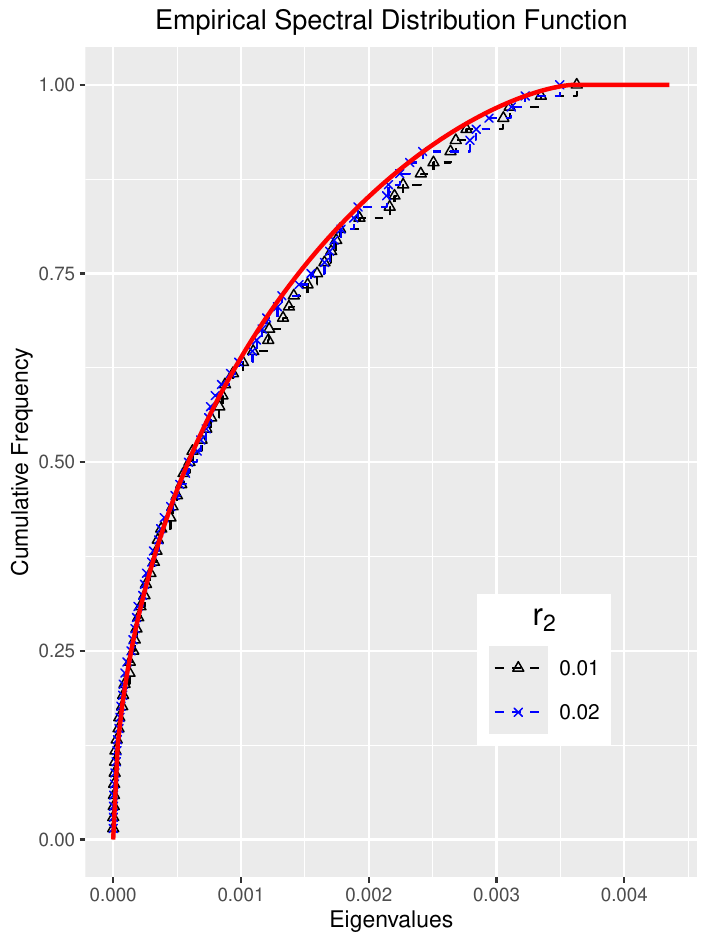}}
	\subfigure[$p=102, k_n=68$]{\includegraphics[width=5cm,height=4cm]{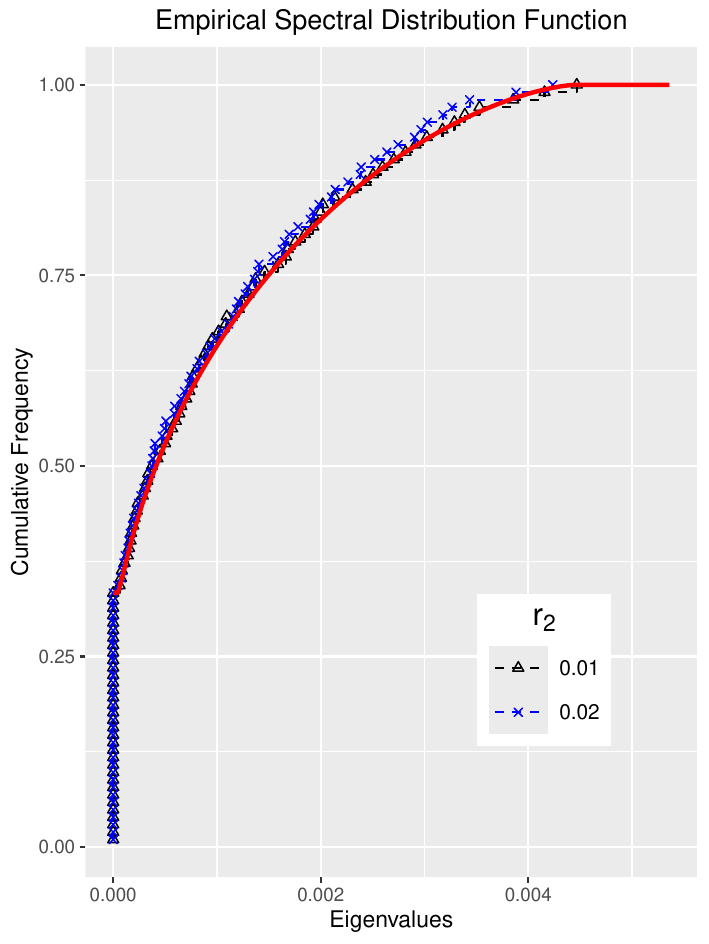}}
	\subfigure[$p=34, k_n=68$]{\includegraphics[width=5cm,height=4cm]{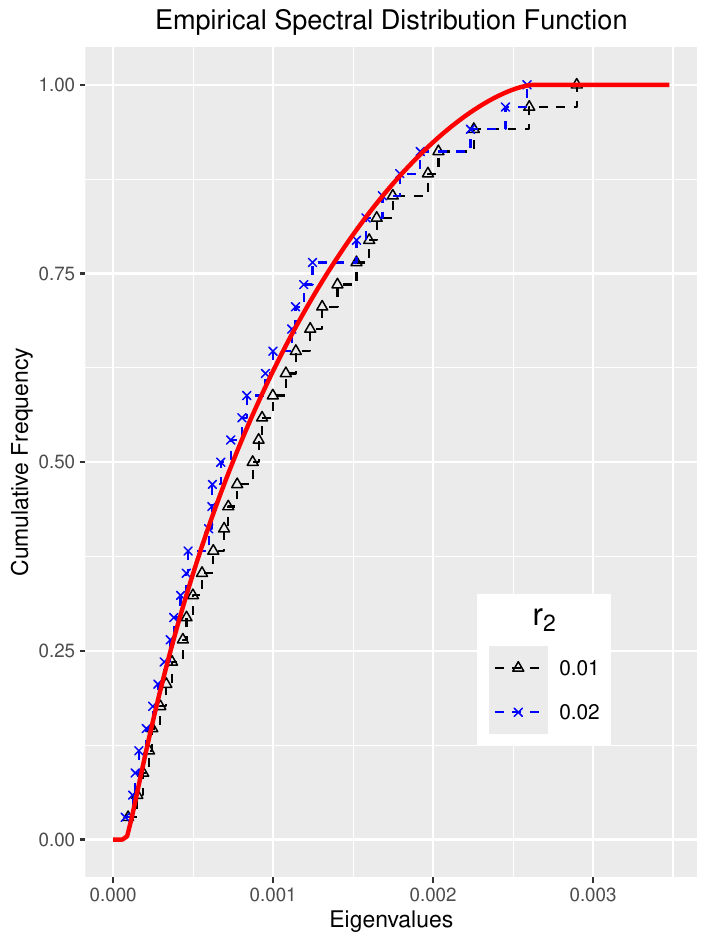}}
	\subfigure[$p=68, k_n=68$]{\includegraphics[width=5cm,height=4cm]{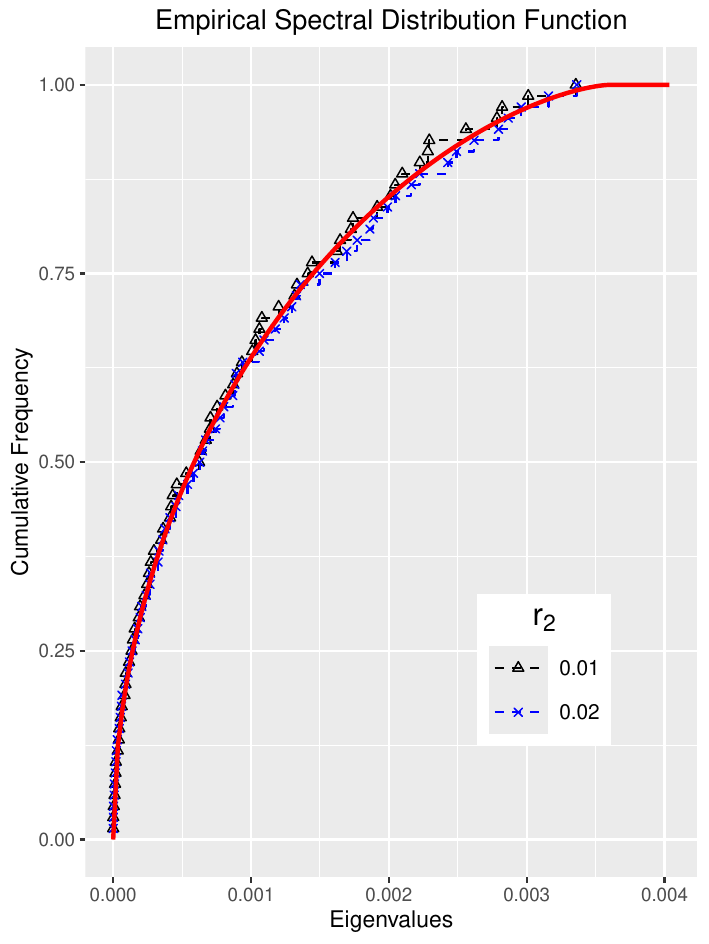}}
	\subfigure[$p=102, k_n=68$]{\includegraphics[width=5cm,height=4cm]{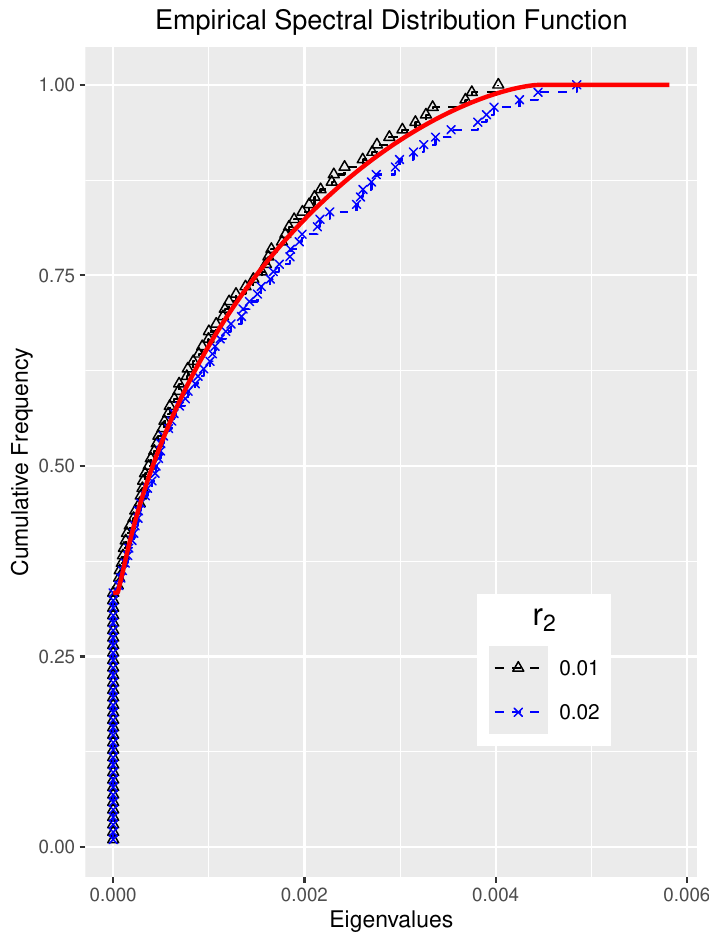}}
	\caption{The empirical spectral distribution functions of realized spot volatility matrix at $t=0$ for the stochastic volatility model \eqref{sto:vol} with $\rho=0.4$ (upper) and  $\rho=0.8$ (lower), with fixed $k_n=68$ and $p=34,68,102$ (Correspondingly, $\bar{p} = 0.5 ,1, 1.5$.). The red solid lines are corresponding theoretical limiting spectral distribution functions, namely the Mar$\check{\text{c}}$enko-Pastur law given by Theorem \ref{spot_esd}.}
	\label{fig_esd2}
\end{figure}

We first present some simulation studies to illustrate the behavior of the empirical spectral distribution of the realized spot volatility matrix estimator, for both the deterministic volatility model \eqref{det:vol} and the stochastic volatility model \eqref{sto:vol}. 
In Figure \ref{fig_esd}--\ref{fig_esd2}, we use red solid lines to represent the cumulative distribution function of Mar$\check{\text{c}}$enko-Pastur law provided by Theorem \ref{spot_esd}, which is also the limiting spectral distribution of the realized spot volatility matrix. Moreover, results for different dimensions $p=34,68,102$ are reported when multiple variational parameters $r_1$ and $r_2$ are used for the deterministic and stochastic volatility models, respectively. Furthermore, different levels of correlation strength $\rho=0.4, 0.8$ is also considered for the stochastic volatility model. 
From the figures, we see that all the empirical spectral distributions are close to corresponding theoretical Mar$\check{\text{c}}$enko-Pastur law, which verifies the conclusions in Theorem \ref{spot_esd}. 

\begin{figure}[!htbp]
	\centering 
	\subfigure[$p=34, k_n=68$]{\includegraphics[width=5cm,height=4cm]{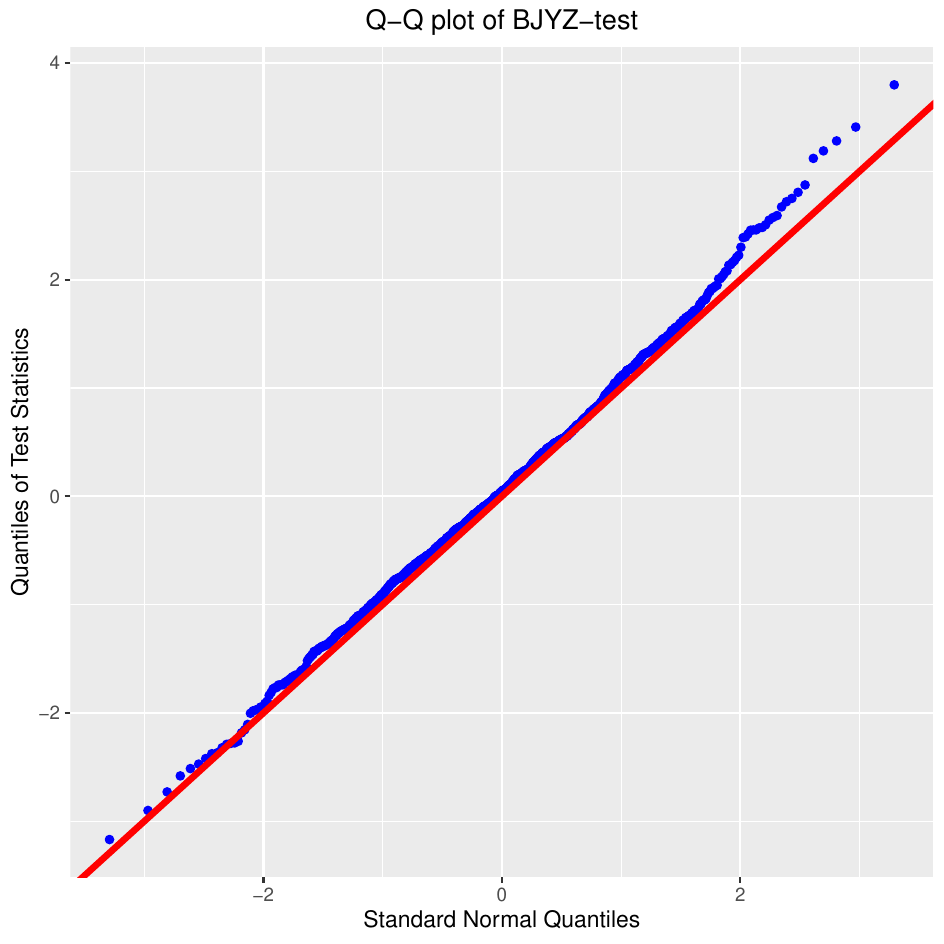}}
	\subfigure[$p=34, k_n=68$]{\includegraphics[width=5cm,height=4cm]{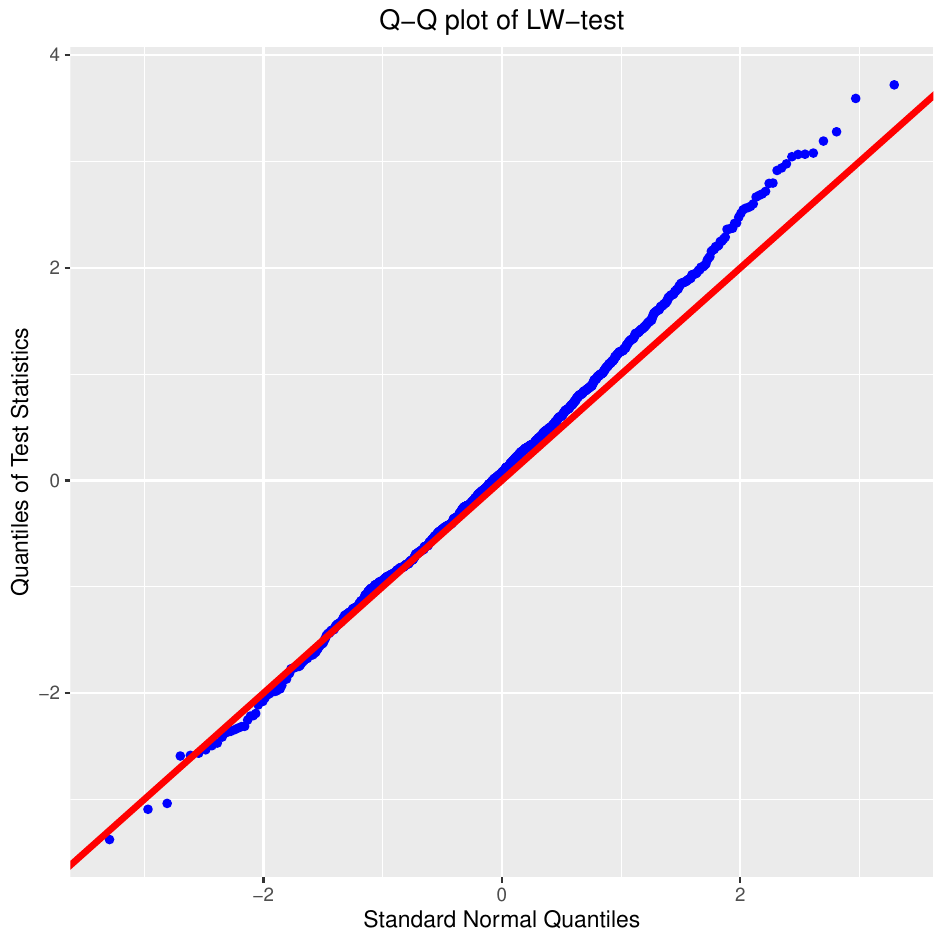}}
	\caption{Q-Q plot for BJYZ-test-spot (left) and LW-test-spot (right) when $\bar{p} = 0.5$ (Namely, $p=34, k_n=68$), based on 1000 repetitions.}
	\label{fig_qq1}
\end{figure}

\begin{figure}[!htbp]
	\centering 
	\subfigure[$p=68, k_n=68$]{\includegraphics[width=5cm,height=4cm]{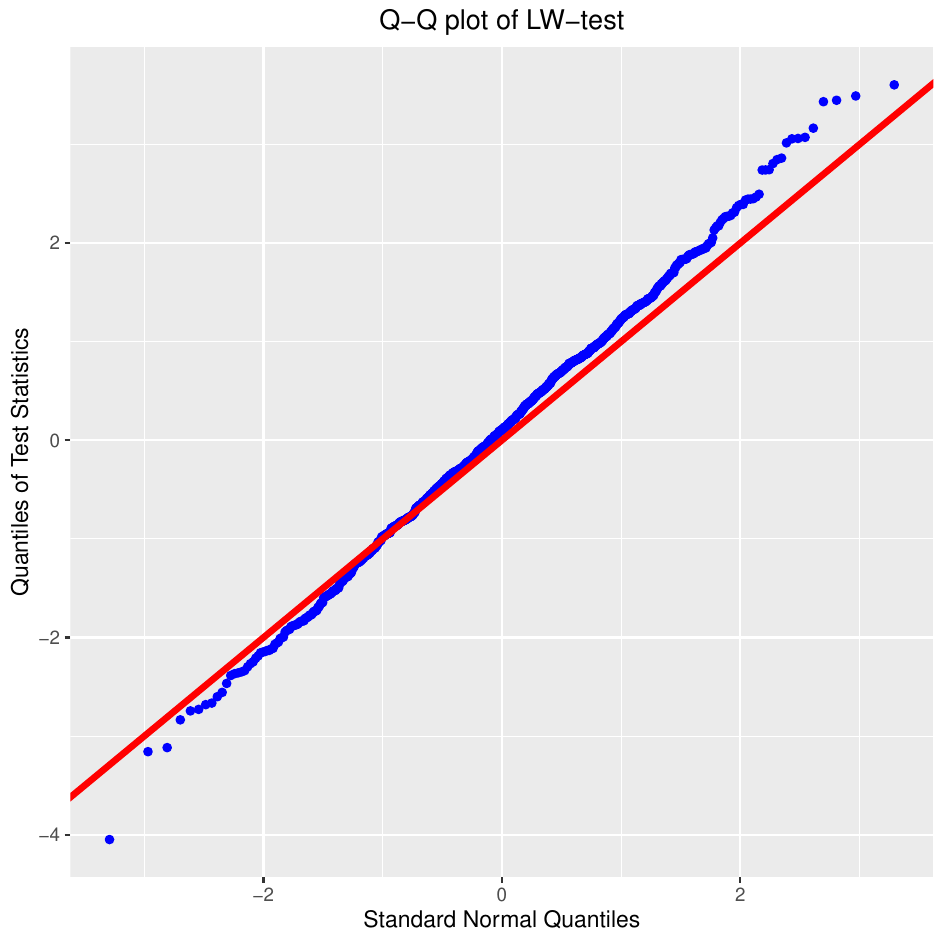}}
	\subfigure[$p=102, k_n=68$]{\includegraphics[width=5cm,height=4cm]{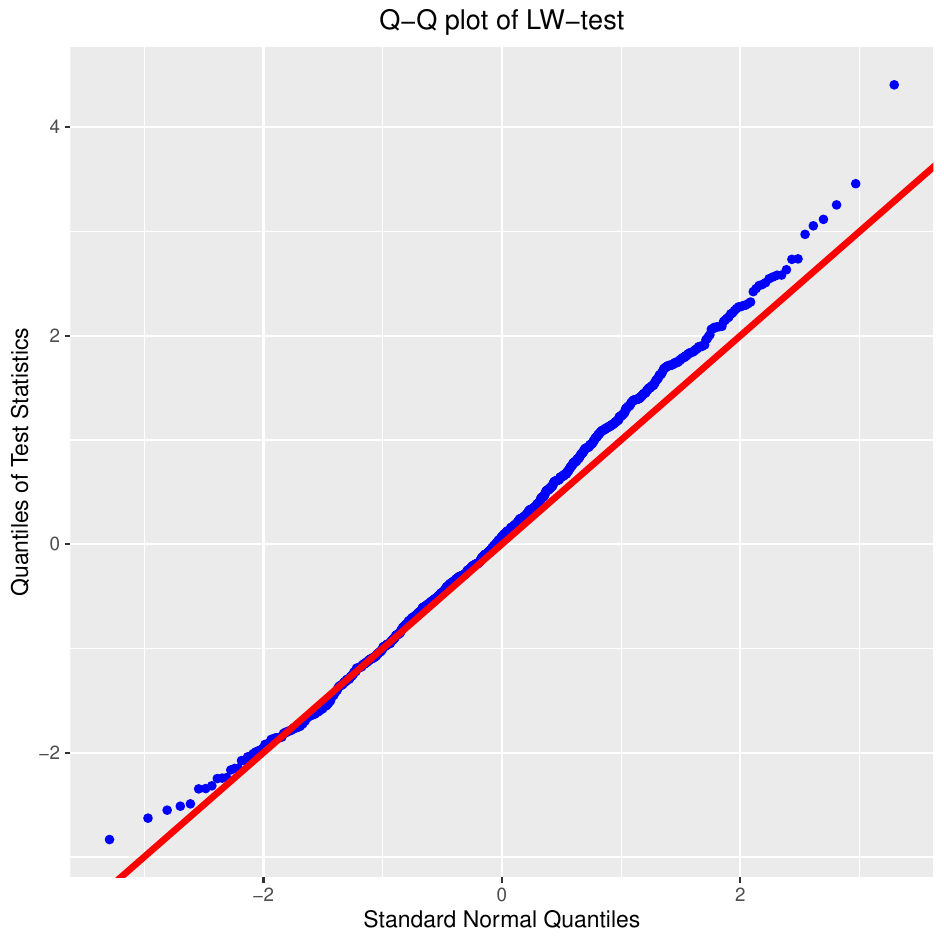}}
	\caption{Q-Q plot for LW-test-spot when $p = 68$ (left) and $p=102$ (right) with fixed $k_n=68$, corresponding to $\bar{p} = 1, 1.5$ respectively, based on 1000 repetitions.}
	\label{fig_qq2}
\end{figure}

\begin{figure}[!htbp]
	\centering 
	\subfigure[$p=34, k_n=68$]{\includegraphics[width=5cm,height=4cm]{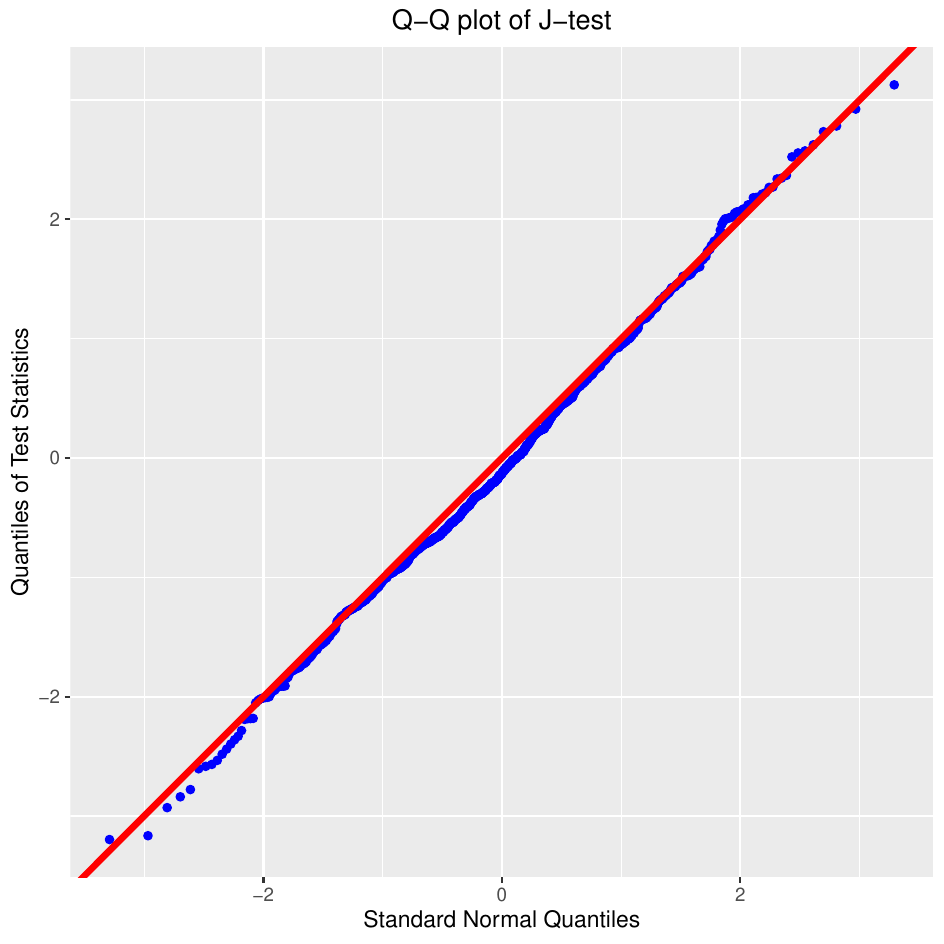}}
	\subfigure[$p=68, k_n=68$]{\includegraphics[width=5cm,height=4cm]{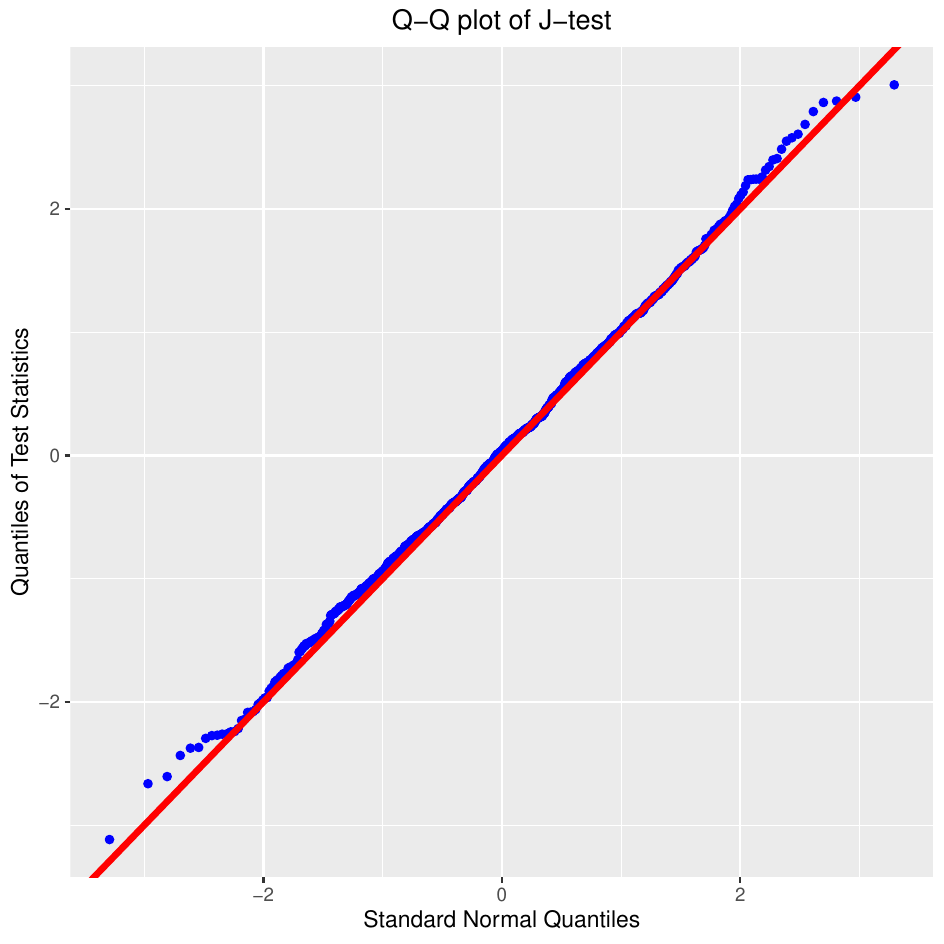}}
	\subfigure[$p=102, k_n=68$]{\includegraphics[width=5cm,height=4cm]{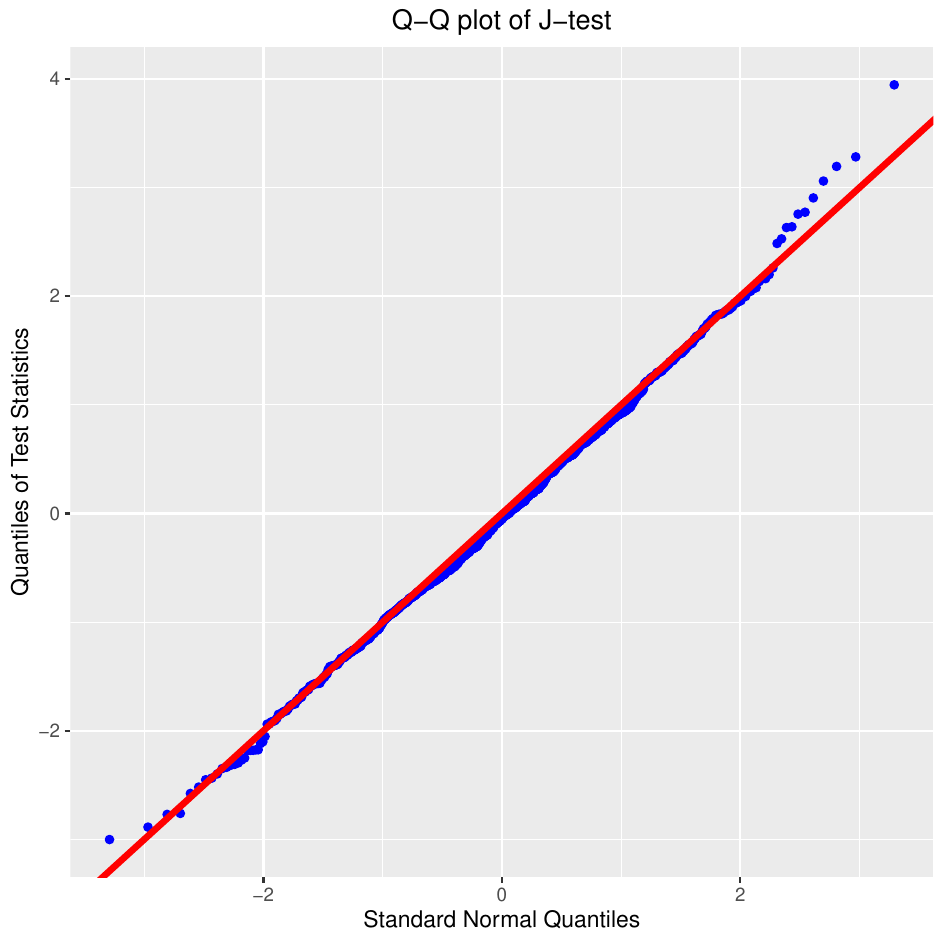}}
	\caption{Q-Q plot for J-test-spot when $p = 34$ (left), $p = 68$ (middle) and $p=102$ (right) with fixed $k_n=68$, corresponding to $\bar{p} = 0.5, 1, 1.5$ respectively, based on 1000 repetitions.}
	\label{fig_qq3}
\end{figure}

Next, we consider the following hypothesis testing problem:
\begin{align}
	H_0: c_{0} = 0.0009 \cdot \mathbb{I}_{p}  \qquad H_1: c_{0} \neq  0.0009\cdot \mathbb{I}_{p},
\end{align}
or equivalently
\begin{align}\label{emp:test}
	H_0: \frac{c_{0}}{0.0009} =  \mathbb{I}_{p}  \qquad H_1: \frac{c_{0}}{0.0009} \neq   \mathbb{I}_{p},
\end{align}
by using the test statistics proposed in Section \ref{sec:app}, and use the deterministic volatility \eqref{det:vol} for illustration. With fixed $r_1 = 0.0008$, we generate a total number of 1000 sample paths and calculate ``BJYZ-test-spot" statistic for $\bar{p}=0.5$ and ``LW-test-spot" statistic for $\bar{p} = 0.5, 1, 1.5$. The Q-Q plots of the resulting estimates are shown in Figure \ref{fig_qq1}--\ref{fig_qq2}. We see that the distributions of the estimates are close to standard normal distribution, which verifies the central limit theorems established in Theorem \ref{thm:testing} and Theorem \ref{thm:testing:lw}. For the hypothesis testing problem \eqref{emp:test}, the empirical sizes of ``BJYZ-test-spot" and ``LW-test-spot" statistics at various significance levels are presented in Table \ref{tab:size}. We observe that when the null hypothesis is true, the sizes of all test statistics are close to corresponding nominal levels. 
Moreover, for $\bar{p} = 0.5$, we see that ``BJYZ-test-spot" generally outperforms ``LW-test-spot", which suggests that the former one should be used if we have prior knowledge of $\bar{p}<1$. 
For both test statistics under different settings of significance level and $\bar{p}$, the size performance is relative worse for larger value of $r_1$. This is to our expectation since as the magnitude of $r_1$ increases, the sample data will gradually deviate from the identically distributed case.
For the evaluation of power performance, we consider the alternative hypothesis of $c_0 = \text{diag}(0.0009\mathbb{I}_{\lfloor sp \rfloor}, 0.0004\mathbb{I}_{p-\lfloor sp \rfloor})$ for $s=0.45, 0.6, 0.75$, namely the diagonal matrix $c_0$ has respectively $s$ of its diagonal elements being 0.0009 whereas the rest $(1-s)$ part are 0.0004. 
Meanwhile, the same data generating scheme as model \eqref{det:vol} is used with $r_1=0, 0.0002, 0.0004$. Table \ref{tab:power} presents the finite sample performances of the test statistics under different settings. The simulation results demonstrate satisfying power performance when $s$ is not larger than 0.6. Besides, we can also see that ``BJYZ-test-spot" generally outperforms ``LW-test-spot" when $\bar{p} = 0.5$, and for both statistics, the power result is relatively worse for larger value of $r_1$. These findings have already been obtained for size performance, and the reasons are same. 
We use the same model and settings for ``J-test-spot" statistic proposed in Section \ref{sec:sphetest} for the sphericity test, and present the Q-Q plots of the estimates in Figure \ref{fig_qq3}, when $\bar{p} = 0.5, 1, 1.5$. We see that the distributions of the estimates are close to standard normal distribution, which verifies the central limit theorems established in Theorem \ref{thm:testing:sph}. As for the sphericity hypothesis testing problem of \eqref{testing:sphe}, the same statistic was considered in \cite{YZC2021} and evaluated with simulation studies under class $\mathcal{C}$ \eqref{vol:classC}, thus we do not repeat the experiment.

\begin{table}[!htbp]
	\caption{Empirical size ($\%$) of BJYZ-test-spot and LW-test-spot based on 1000 repetitions.}
	\centering
	\scriptsize
	\begin{tabular}{|c|c|c|c|c|c|c|c|c|c|c|c|}
		\toprule
		\multirow{2}{*}{Significance level}	& \multirow{2}{*}{Test statistics} &
		\multicolumn{3}{c|}{$\text{r}_1 = 0$}& \multicolumn{3}{c|}{$\text{r}_1 = 0.0004$}& \multicolumn{3}{c|}{$\text{r}_1 = 0.0008$}\\
		\cline{3-11}& 
		& $\bar{p}=0.5$ &  $\bar{p}=1$  &  $\bar{p}=1.5$ & $\bar{p}=0.5$ &  $\bar{p}=1$  &  $\bar{p}=1.5$ & $\bar{p}=0.5$ &  $\bar{p}=1$  &  $\bar{p}=1.5$ \\
		\hline
		\multirow{2}{*}{$10\%$}	& 
		BJYZ-test-spot &10.0 &  ------& ------& 10.3 &  ------& ------& 11.2 &  ------& ------\\
		\cline{2-11}&
		LW-test-spot  & 9.7&9.3& 10.9& 11.6&11.1&11.6 &14.5&14.8& 14.8\\
		\hline
		\multirow{2}{*}{$5\%$}	& 
		BJYZ-test-spot &  5.5& ------ &------ & 5.8 & ------ &------ & 5.8 &  ------&------ \\
		\cline{2-11}&
		LW-test-spot  &5.1&5.4& 4.8&6.6&5.6& 6.9 &8.7&9&8.2 \\
		\hline
		\multirow{2}{*}{$1\%$}	& 
		BJYZ-test-spot & 0.6 &------  & ------& 1.6&  ------&------ & 2.1&------  &------ \\
		\cline{2-11}&
		LW-test-spot  &1.4&1.3& 0.8&1.6&1.3&1.5 &2.6&2.4& 2.2\\
		\bottomrule
	\end{tabular}\label{tab:size}
\end{table}

\begin{table}[!htbp]
	\caption{Empirical power ($\%$) of BJYZ-test-spot and LW-test-spot based on 1000 repetitions.}
	\centering
	\scriptsize
	\begin{tabular}{|c|c|c|c|c|c|c|c|c|c|c|c|}
		\toprule
		\multicolumn{11}{|c|}{$s=0.45$} \\
		\hline
		\multirow{2}{*}{Significance level}	& \multirow{2}{*}{Test statistics} &
		\multicolumn{3}{c|}{$\text{r}_1 = 0$}& \multicolumn{3}{c|}{$\text{r}_1 = 0.0002$}& \multicolumn{3}{c|}{$\text{r}_1 = 0.0004$}\\
		\cline{3-11}& 
		& $\bar{p}=0.5$ &  $\bar{p}=1$  &  $\bar{p}=1.5$ & $\bar{p}=0.5$ &  $\bar{p}=1$  &  $\bar{p}=1.5$ & $\bar{p}=0.5$ &  $\bar{p}=1$  &  $\bar{p}=1.5$ \\
		\hline
		\multirow{2}{*}{$10\%$}	& 
		BJYZ-test-spot & 100 &  ------& ------& 100 &------  & ------&  100 & ------ & ------\\
		\cline{2-11}&
		LW-test-spot  &100&100&100 &100&100&100 &100&100&100 \\
		\hline
		\multirow{2}{*}{$5\%$}	& 
		BJYZ-test-spot &  100 & ------ & ------&100  &  ------& ------& 100 & ------ & ------\\
		\cline{2-11}&
		LW-test-spot  &100&100& 100&100&100&100 &100&100&100 \\
		\hline
		\multirow{2}{*}{$1\%$}	& 
		BJYZ-test-spot & 100 &------  &------ & 100 & ------ & ------& 100 &  ------&------ \\
		\cline{2-11}&
		LW-test-spot  &100&100&100 &100&100& 100&100&100& 100\\
		\toprule
		\multicolumn{11}{|c|}{$s=0.6$} \\
		\hline
		\multirow{2}{*}{Significance level}	& \multirow{2}{*}{Test statistics} &
		\multicolumn{3}{c|}{$\text{r}_1 = 0$}& \multicolumn{3}{c|}{$\text{r}_1 = 0.0002$}& \multicolumn{3}{c|}{$\text{r}_1 = 0.0004$}\\
		\cline{3-11}& 
		& $\bar{p}=0.5$ &  $\bar{p}=1$  &  $\bar{p}=1.5$ & $\bar{p}=0.5$ &  $\bar{p}=1$  &  $\bar{p}=1.5$ & $\bar{p}=0.5$ &  $\bar{p}=1$  &  $\bar{p}=1.5$ \\
		\hline
		\multirow{2}{*}{$10\%$}	& 
		BJYZ-test-spot & 100 &------  & ------& 99.9 &  ------& ------& 100 &------  & ------\\
		\cline{2-11}&
		LW-test-spot  &100 &100& 100&99.8&100&99.9 &99.9&100& 99.7\\
		\hline
		\multirow{2}{*}{$5\%$}	& 
		BJYZ-test-spot & 100 &  ------&------ & 99.9 &  ------&------ &  100& ------ & ------\\
		\cline{2-11}&
		LW-test-spot  & 100 &99.9& 99.9&99.7&99.9&99.7 &99.6&99.7&99.3 \\
		\hline
		\multirow{2}{*}{$1\%$}	& 
		BJYZ-test-spot & 99.8 & ------ & ------&  99.5&  ------& ------& 99.3 &  ------&------ \\
		\cline{2-11}&
		LW-test-spot  & 98.8 &98.9&98.3 &97.3&98.9&96.8 &95.2&95.8&94.8 \\
		\toprule
		\multicolumn{11}{|c|}{$s=0.75$} \\
		\hline
		\multirow{2}{*}{Significance level}	& \multirow{2}{*}{Test statistics} &
		\multicolumn{3}{c|}{$\text{r}_1 = 0$}& \multicolumn{3}{c|}{$\text{r}_1 = 0.0002$}& \multicolumn{3}{c|}{$\text{r}_1 = 0.0004$}\\
		\cline{3-11}& 
		& $\bar{p}=0.5$ &  $\bar{p}=1$  &  $\bar{p}=1.5$ & $\bar{p}=0.5$ &  $\bar{p}=1$  &  $\bar{p}=1.5$ & $\bar{p}=0.5$ &  $\bar{p}=1$  &  $\bar{p}=1.5$ \\
		\hline
		\multirow{2}{*}{$10\%$}	& 
		BJYZ-test-spot & 97.4 &------  & ------& 96.4 &  ------& ------& 93.6 &  ------& ------\\
		\cline{2-11}&
		LW-test-spot  &89.5&85.3& 87.1&84.5&81.5&81.6 &83.2&79.1&79.9 \\
		\hline
		\multirow{2}{*}{$5\%$}	& 
		BJYZ-test-spot & 94.0 &------  & ------& 92.3 &  ------& ------& 88.6 & ------ & ------\\
		\cline{2-11}&
		LW-test-spot  &79.6&74.3&77.7 &73.4& 71.3&69.8 &73.5&65.9& 68.9\\
		\hline
		\multirow{2}{*}{$1\%$}	& 
		BJYZ-test-spot & 83.8 &------  & ------& 79.9 & ------ &------ & 74.5 & ------ &------ \\
		\cline{2-11}&
		LW-test-spot  & 52.7 &45.8&49.9 &48.6& 43.6&42.6 &44.7&37.1& 40.5\\
		\bottomrule
	\end{tabular}\label{tab:power}
\end{table}

\section{Conclusion and future work}\label{sec:conc}
In this article, we conduct spectral analysis for the realized spot volatility matrix estimator by using the high-dimensional and high-frequency data, including its limiting spectral distribution and linear spectral statistics. Our theoretical results demonstrate that, under some general conditions, the limiting spectral distributions of the spot volatility matrix and the realized spot volatility matrix estimator are governed by a regular Mar$\check{\text{c}}$enko-Pastur equation. Furthermore, the central limit theorem for linear spectral statistics of the realized spot covariance matrix estimator is established. The result is similar to the one for independent and identically distributed sample in classical large-dimensional random matrix theory. As applications, we perform identity and sphericity tests for the spot volatility matrix and propose some test statistics by using our established theory. Then, the theoretical conclusions and the finite sample performance of the proposed statistics are justified by our simulation studies.

We note that the spectral theory developed for the spot volatility matrix in this paper could potentially be used to solve some other problems. For example, we may consider to test if $c_{t_1} = c_{t_2}$ for $0 \leq t_1< t_2\leq 1$.
Another problem is that, in \cite{ZL2011} and \cite{YZC2021}, the spot volatility matrix is restricted to satisfy a special structure of class $\mathcal{C}$ as in \eqref{vol:classC}, while such a structure remains to be testified by using the real data. From a high-frequency perspective, we may extend the analysis to the scenario with the presence of both jumps and market microstructure noise. \red{In this paper, we do not consider the jumps for simplicity. Nevertheless, as studied by Jacod and Protter (2012) (Theorem 13.3.3), the spot volatility estimator remains consistent when the jumps are not too frequent. For a general L$\acute{\text{e}}$vy jump part, a truncated version for the estimator is available, we refer to Jacod and Protter (2012) for more details.}
Those mentioned problems are out the scope of this paper and will be considered in our future works. 
%
%
%
%
	
	\hspace{-0.25in}
	\bibliographystyle{model2-names}
	\bibliography{liu}
\end{document}